\documentclass[pdflatex,sn-mathphys-num]{sn-jnl}
\usepackage{fullpage}
\usepackage{hyperref}
\usepackage{amssymb}
\usepackage{amsmath}
\usepackage{mdframed}
\usepackage{subfig}
\usepackage{cancel}
\usepackage{enumerate}
\usepackage{color}
\usepackage{mathtools}
\usepackage{pdfpages}
\usepackage{mathrsfs}
\usepackage{enumitem}
\usepackage{mathtools}
\usepackage{graphicx}
\usepackage{mdframed}
\usepackage{amsfonts}
\usepackage{cancel}
\usepackage{placeins}
\usepackage{amsthm}
\usepackage{tikz,epsfig}
\usetikzlibrary{shapes,arrows,automata}
\usepackage{xcolor}
\usepackage[normalem]{ulem}
\usetikzlibrary{arrows.meta}
\numberwithin{equation}{section}
\numberwithin{figure}{section}

\definecolor{goldenpoppy}{rgb}{0.99, 0.76, 0.0}
\definecolor{richblack}{rgb}{0.06, 0.05, 0.03}
\definecolor{cadmiumred}{rgb}{0.89, 0.0, 0.13}
\definecolor{fuchsia}{rgb}{0.3, 0.0, 0.3}
\definecolor{green(ncs)}{rgb}{0.0, 0.52, 0.32}
\tikzstyle{species_T} = [circle,radius=0.1cm, text centered, draw=black, fill=goldenpoppy]
\tikzstyle{species_R} = [circle,radius=0.1cm, text centered, draw=black, fill=white]
\tikzstyle{species_C} = [circle,radius=0.1cm, text centered, draw=black, fill=fuchsia]
\tikzstyle{dots} = [circle,radius=0.1cm,text centered]
\tikzstyle{arrowA} = [-{Latex[length=2mm]},white,dashed]
\tikzstyle{arrowB} = [-{Latex[length=2mm]},green(ncs)]
\tikzstyle{arrowC} = [-{Latex[length=2mm]},cadmiumred]
\tikzstyle{inhibit} = [thick,-|,black!100]
\tikzstyle{loosely dashed}= [dash pattern=on 3pt off 6pt]

\newtheorem{definition}{Definition}[section]
\newtheorem{theorem}{Theorem}[section]

\newtheorem{example}{Example}[section]
\newtheorem{lemma}{Lemma}[section]
\newtheorem{corollary}{Corollary}[section]

\let\Item\item

\mdfdefinestyle{theoremstyle}{%
frametitlerule=true,%
innertopmargin=\topskip,
}
\mdtheorem[style=theoremstyle]{algorithm}{Algorithm}[section]

\begin{document}

\title{\huge From Algebraic Curves to Chemical Limit Cycles}

\author[1,2]{\fnm{Tomislav} \sur{Plesa}}\email{tp525@cam.ac.uk}

\author[3]{\fnm{Radek} \sur{Erban}} \email{radek.erban@maths.ox.ac.uk}

\author[4]{\fnm{Armengol} \sur{Gasull} \email{armengol.gasull@uab.cat}} 

\author[5]{\fnm{Gheorghe} \sur{Craciun} \email{craciun@math.wisc.edu}}

\affil[1]{\orgdiv{Department of Applied Mathematics and Theoretical Physics},
\orgname{University of Cambridge},
\orgaddress{\city{Cambridge}, \postcode{CB3 0WA}, \country{United Kingdom}}}

\affil[2]{\orgdiv{School of Mathematical Sciences},
\orgname{Lancaster University},
\orgaddress{\city{Lancaster}, \postcode{LA1 4YF}, \country{United Kingdom}}}

\affil[3]{\orgdiv{Mathematical Institute}, \orgname{University of Oxford}, \orgaddress{\street{Radcliffe Observatory Quarter,
Woodstock Road}, \city{Oxford}, \postcode{OX2 6GG}, \country{United Kingdom}}}

\affil[4]{\orgdiv{Departament de Matemàtiques,
Facultat de Ciències},
\orgname{Universitat Autònoma de Barcelona},
\orgaddress{\postcode{08193} \city{Cerdanyola del Vallés, Barcelona}, \country{Spain}}}

\affil[5]{\orgdiv{Department of Mathematics and
Department of Biomolecular Chemistry},
\orgname{University of Wisconsin-Madison},
\orgaddress{\street 405 Van Vleck Hall, 480 Lincoln Dr},
\city{Madison}, \postcode{WI 53706-1388}, \country{USA}}

\abstract{Algebraic limit cycles in planar chemical dynamical systems 
are investigated through their connection with real algebraic curves. 
General constructions are established for realizing 
algebraic ovals as hyperbolic algebraic limit cycles 
of chemical dynamical systems. This framework is 
used to construct a cubic chemical system 
with two algebraic limit cycles, 
a quartic system with four algebraic limit cycles, 
and a family of systems in which the number of algebraic limit 
cycles grows quadratically with the system degree.
The framework is then specialized to chemical dynamical systems 
admitting reversible reaction-network realizations. 
In particular, the framework yields 
reversible cubic chemical systems with an algebraic limit cycle, 
a reversible quartic system with two algebraic limit cycles, 
and a family of reversible systems with 
quadratically many algebraic limit cycles relative to their degree.
These results improve upon existing lower bounds for the maximal number of algebraic limit cycles in chemical dynamical systems of a given degree.
Moreover, this work identifies special subclasses of polynomials and Harnack M-curves, 
and establishes connections between the topology of real algebraic curves and 
oscillatory dynamics in chemical reaction networks.}


\maketitle


\section{Introduction}
\label{sec:1}

Limit cycles are among the fundamental objects in the qualitative theory of planar dynamical systems, and determining how many limit cycles can occur in polynomial systems is a classical and challenging problem~\cite{Hilbert,Ilyashenko,Perko}. A particularly tractable class is formed by algebraic limit cycles, namely periodic orbits contained in real algebraic curves~\cite{Llibre,Llibre2,Gasull_Santana}. In chemical dynamical systems, this question acquires an additional structure: the polynomial vector field is constrained by the requirement that it admit a realization as a chemical reaction network~\cite{QCM}. The purpose of this paper is to investigate how chemical constraints affect the existence and number of algebraic limit cycles, both in general chemical dynamical systems and in the more restrictive subclass of systems admitting reversible reaction-network realizations~\cite{Feinbergbook,Janos}.

The study of algebraic limit cycles is closely related 
to the topology of real algebraic curves. 
This connection was investigated by Winkel~\cite{Winkel},
who showed that the ovals of a regular (non-singular) algebraic curve of degree $n$ can be realized as attracting limit cycles of a polynomial vector field of degree $2n-1$. This construction was subsequently improved by Christopher~\cite{Christopher}, who obtained polynomial vector fields of degree $n$ having precisely the ovals  of the prescribed $n$-degree algebraic curve as limit cycles. 
To maximize the number of algebraic limit cycles obtained through such constructions, 
it is therefore natural to consider algebraic curves with the 
largest possible number of ovals. A classical result of Harnack~\cite{Harnack} establishes an upper bound $M(n)$ on the number of 
bounded connected components of a regular 
real plane algebraic curve of a degree $n$, namely
\begin{align}
M(n) & = 
\begin{cases}
\frac{(n-1)(n-2)}{2} + 1, & \text{for every even } n \ge 2, \\
 \frac{(n-1)(n-2)}{2}, & \text{for every odd } n \ge 1.
\end{cases}
\label{eq:M_n}
\end{align}
This bound is attained by the so-called Harnack M-curves.

More broadly, the connection between algebraic curves
and limit cycles of planar polynomial dynamical systems is
reflected in Hilbert's 16th problem~\cite{Hilbert}, 
which consists of two parts: the first concerns the 
possible arrangements of ovals of real algebraic curves of a given degree, 
while the second part concerns the maximum number and possible 
arrangements of limit cycles of planar polynomial dynamical systems of a given degree. 
Algebraic limit cycles provide a direct connection between the two parts
of Hilbert's 16th problem, being both ovals of algebraic curves
and periodic orbits of polynomial dynamical systems.

Both parts of Hilbert's 16th problem remain unsolved in general,
as does the restriction of its second part to algebraic limit cycles.
In particular, let $H_a(n)$ denote the supremum of the number 
of algebraic limit cycles in the real plane $\mathbb{R}^2$ 
among planar polynomial dynamical systems of degree at most $n$.
It was shown in~\cite{Llibre,Llibre2} that, for planar polynomial dynamical systems
of degree $n$ whose invariant algebraic curves satisfy suitable genericity conditions,
the maximum number of algebraic limit cycles is precisely the Harnack number 
$M(n)$ given by~(\ref{eq:M_n}), implying $H_a(n) \ge M(n)$.
However, more than $M(n)$ algebraic limit cycles may exist
when the genericity conditions are violated.
In particular, cubic polynomial systems with 
two algebraic limit cycles have been constructed~\cite{Llibre_Zhao},
yielding $H_a(3) \ge 2 > M(3) = 1$. 

More recently, the corresponding problem for chemical reaction networks was investigated in~\cite{CEpaper}, where lower bounds on the number of algebraic and non-algebraic limit cycles in several classes of planar chemical systems, including those with reversible reaction networks, were presented. Constructions of chemical systems with multiple prescribed algebraic limit cycles were obtained, and the relationship between the degree of the underlying algebraic curve and that of the corresponding chemical dynamical system was investigated. In particular, let $H_a^c(n)$ and $H_{a,r}^{c}$ 
denote suprema of the number of algebraic limit cycles 
in the positive quadrant $\mathbb{R}_{>}^2$ 
among $n$-degree planar chemical dynamical systems
and those admitting reversible reaction-network realizations, respectively.
Clearly, $H_{a,r}^{c}(n) \le H_{a}^{c}(n) \le H_{a}(n)$.
In Table~\ref{tab:1}, we summarize the lower bounds for the Hilbert numbers 
$H_a^c(n)$ and $H_{a,r}^{c}$ presented in~\cite{CEpaper}, 
together with the known lower bounds for $H_a(n)$.

\begin{table}[h!]
\begin{tabular}{|c|c|c|c|}
\hline
\raise -1.4mm \hbox{\rule{0pt}{4.6mm}}
$\quad\;\;\,$ {\color{black} $n$} $\quad\quad$
&
$\quad\quad$ {\color{black} $H_{a,r}^{c}(n)$} $\quad\quad$
&
$\quad\quad${\color{black} $H_{a}^c(n)$} $\quad\quad$
&
$\quad\quad$ {\color{black} $H_{a}(n)$} $\quad\quad$
\\
\hline
\hline
\raise -1.4mm \hbox{\rule{0pt}{4.6mm}}
{\color{black} 2} & = 0 & $\ge$ 1 & $\ge$ 1 \\
\hline
\raise -1.4mm \hbox{\rule{0pt}{4.6mm}}
{\color{black} 3} & $\ge$ 0 & $\ge$ 1 & $\ge$ 2 \\
\hline
\raise -1.4mm \hbox{\rule{0pt}{4.6mm}}
{\color{black} 4} & $\ge$ 1 & $\ge$ 3 & $\ge$ 4 \\
\hline
\hline
\hbox{\rule{0pt}{4.6mm}}
{\color{black} $n \ge 5$}
&
$\ge \left\lfloor \frac{n-2}{4} \right\rfloor$
&
$\ge M(n-2)$
&
$\ge M(n)$
\\
\hline
\end{tabular}
\vskip 2mm
\caption{{\it Lower bounds for the Hilbert numbers 
$H_{a,r}^{c}(n)$ and $H_{a}^{c}(n)$ presented in~\emph{\cite{CEpaper}}.
Also presented in~\emph{\cite{CEpaper}} are results 
for other chemical Hilbert numbers; 
see \emph{Section~\ref{sec:6}} for a discussion.}}
\label{tab:1}
\end{table}

\begin{table}[h!]
\begin{tabular}{|c|c|c|c|}
\hline
\raise -1.4mm \hbox{\rule{0pt}{4.6mm}}
$\quad\;\;\,$ {\color{black} $n$} $\quad\quad$
&
$\quad\quad$ {\color{black} $H_{a,r}^{c}(n)$} $\quad\quad$
&
$\quad\quad${\color{black} $H_{a}^c(n)$} $\quad\quad$
&
$\quad\quad$ {\color{black} $H_{a}(n)$} $\quad\quad$
\\
\hline
\hline
\raise -1.4mm \hbox{\rule{0pt}{4.6mm}}
{\color{black} 2} & = 0 & $\ge$ 1 & $\ge$ 1 \\
\hline
\raise -1.4mm \hbox{\rule{0pt}{4.6mm}}
{\color{black} 3} & {\color{red} $\ge$ 1} & {\color{red} $\ge$ 2} & $\ge$ 2 \\
\hline
\raise -1.4mm \hbox{\rule{0pt}{4.6mm}}
{\color{black} 4} & {\color{red} $\ge$ 2} &{\color{red} $\ge$ 4} & $\ge$ 4 \\
\hline
\hline
\hbox{\rule{0pt}{4.6mm}}
{\color{black} $n \ge 5$}
&
{\color{red} $\ge M(n-3)$}
&
{\color{red}  $\ge M(n) - \frac{1 + (-1)^n}{2}$}
&
$\ge M(n)$
\\
\hline
\end{tabular}
\vskip 2mm
\caption{{\it Lower bounds for the Hilbert numbers 
$H_{a,r}^{c}(n)$ and $H_{a}^{c}(n)$ established in this paper.
The bounds $H_{a}^c(3) \ge 2$, $H_{a}^c(4) \ge 4$, 
$H_{a}^c(n) \ge M(n) - (1 + (-1)^n)/2$ for $n \ge 5$ are established 
in \emph{Corollaries~\ref{corollary:cubic_2}} and~\emph{\ref{corollary:quartic_4}},
and \emph{Theorem~\ref{theorem:semi_chemical_improved}}, respectively.
The bounds $H_{a,r}^c(3) \ge 1$, $H_{a,r}^c(4) \ge 2$, 
$H_{a,r}^c(n) \ge M(n-3)$ for $n \ge 5$ are established 
in \emph{Corollaries~\ref{corollary:cubic_R}} and~\emph{\ref{corollary:quartic_2_R}},
and \emph{Theorem~\ref{theorem:buffered_n}}, respectively.}}
\label{tab:2}
\end{table}

\begin{figure}[!htbp]
\vskip 0.0cm
\leftline{\hskip 2.0cm (a) CDS~(\ref{eq:cubic_2}) 
\hskip  4.6cm (b) CDS~(\ref{eq:quartic_4})}
\vskip 0.2cm
\centerline{
\hskip -1.0cm
\includegraphics[width=0.35\columnwidth]{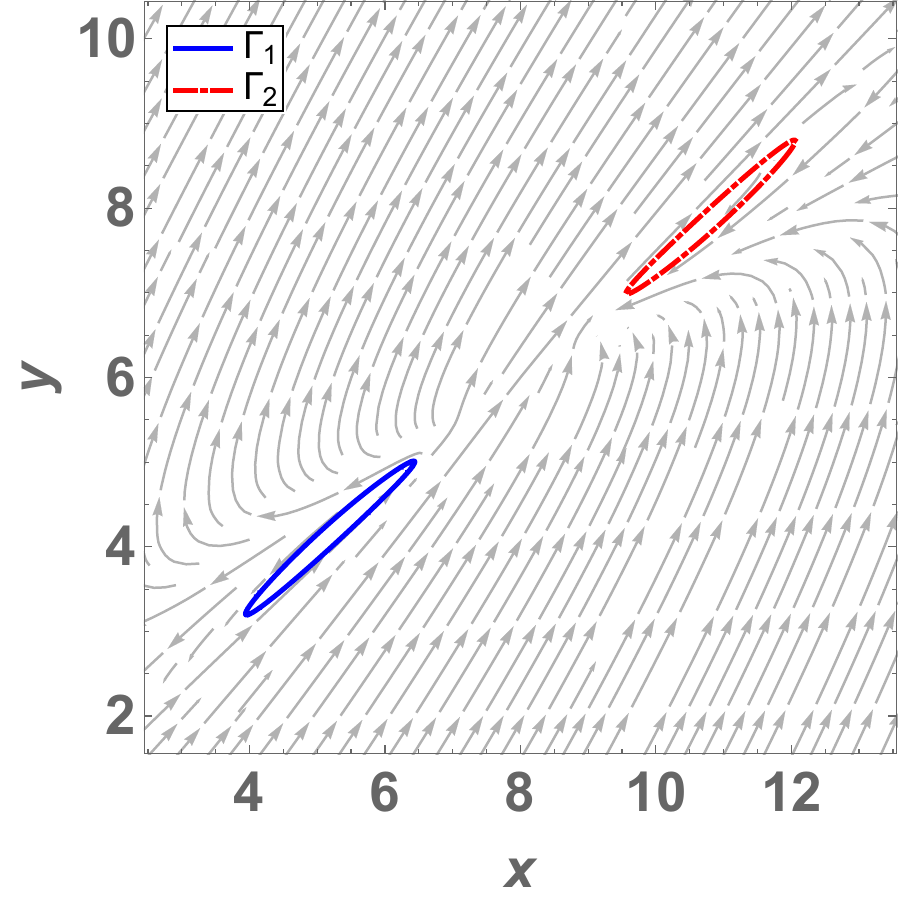}
\hskip 1.2cm
\includegraphics[width=0.35\columnwidth]{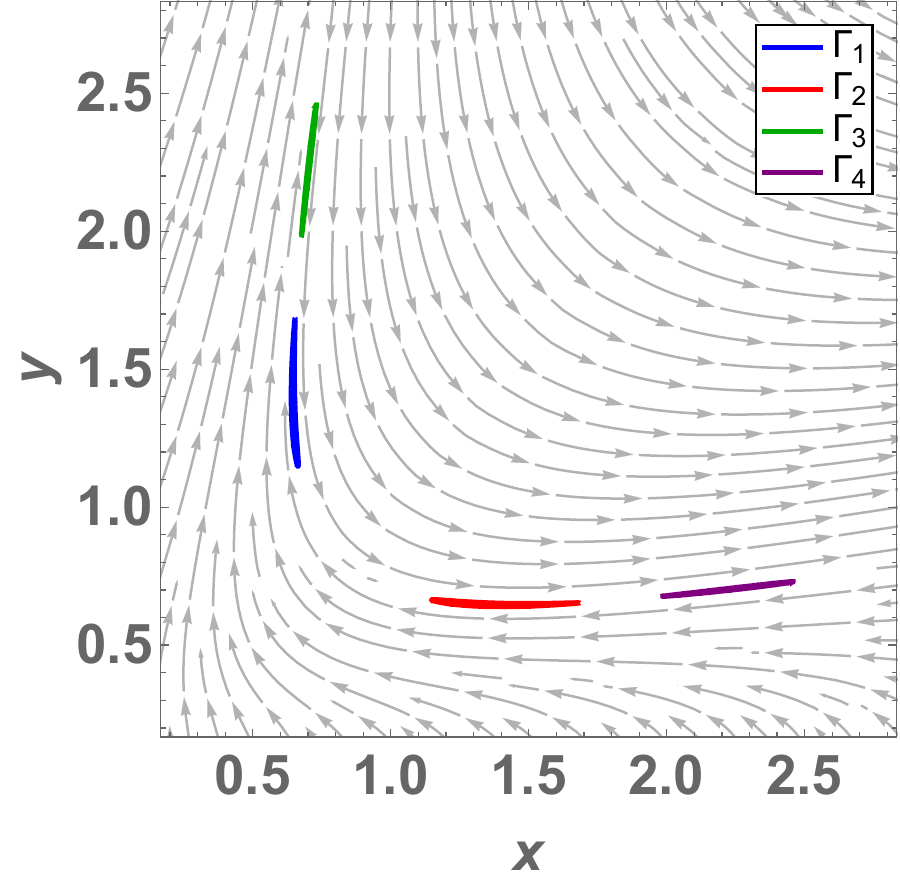}
}
\leftline{\hskip 2.0cm (c) CDS~(\ref{eq:cubic_R})
\hskip  4.6cm (d) CDS~(\ref{eq:quartic_2_R})}
\vskip 0.2cm
\centerline{
\hskip -1.0cm
\includegraphics[width=0.35\columnwidth]{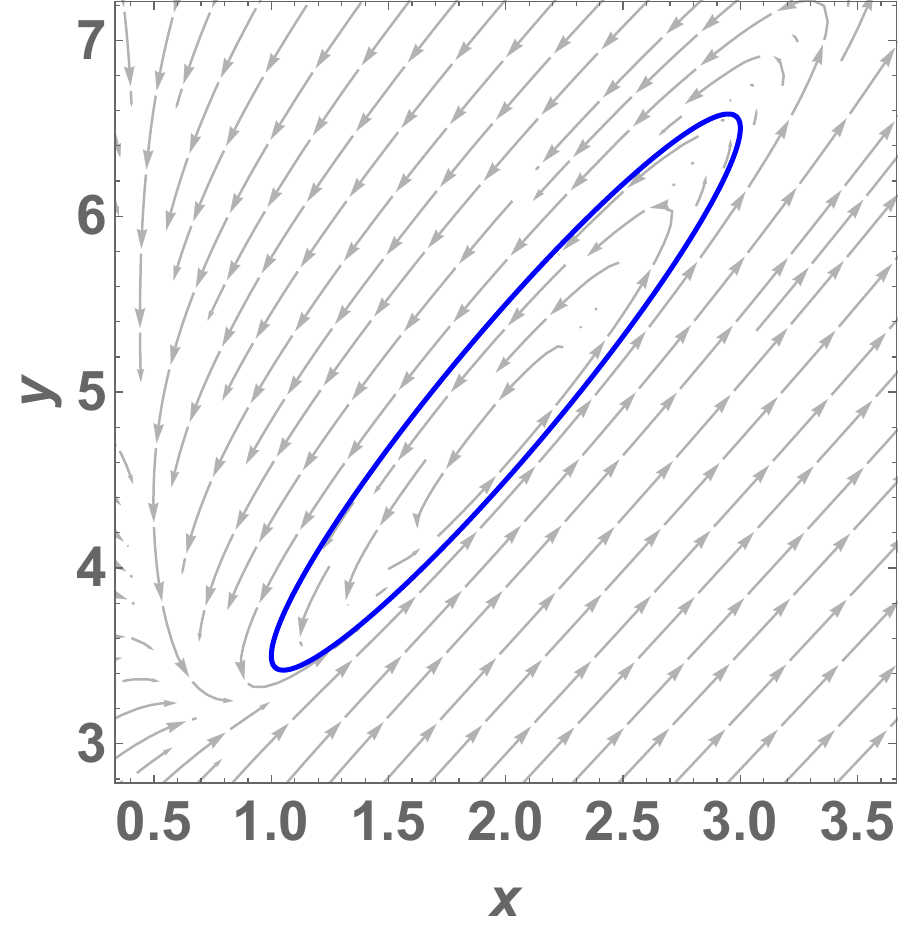}
\hskip 1.2cm
\includegraphics[width=0.35\columnwidth]{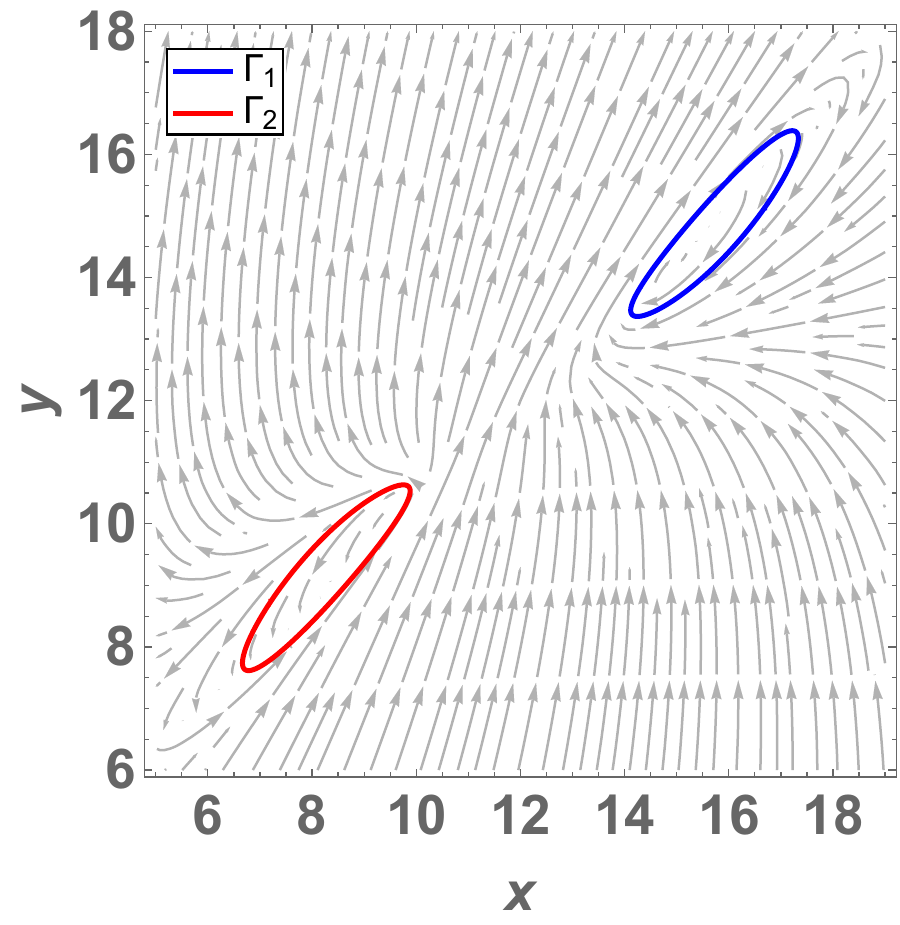}
}
\vskip 0.0cm
\caption{\it Panels \emph{(a)} and \emph{(b)} display phase-plane diagrams
of the cubic chemical dynamical system~$(\ref{eq:cubic_2})$,
with two algebraic limit cycles (one stable, shown by a solid blue curve, 
and one unstable, shown by a dashed-dotted red curve), and 
of the quartic system~$(\ref{eq:quartic_4})$, with $\varepsilon = 5/100$
and four stable algebraic limit cycles, respectively.
The four limit cycles in panel \emph{(b)} appear flattened at this scale; 
see Supplementary Information for zoomed-in plots confirming they are ovals. 
Panels \emph{(c)} and \emph{(d)} display diagrams 
for the reversible cubic chemical dynamical system~$(\ref{eq:cubic_R})$,
with $\mu = 1/100$ and one stable algebraic limit cycle,
and the reversible quartic system~$(\ref{eq:quartic_2_R})$,
with two stable algebraic limit cycles, respectively.} 
\label{fig:1}
\end{figure}

In this paper, we extend the results from~\cite{Winkel,Christopher}
and develop general constructions for realizing 
algebraic ovals as algebraic cycles of chemical dynamical systems, 
including conditions ensuring their hyperbolicity and criteria
for determining their stability and excluding additional cycles.
We use this framework to connect the topology of algebraic curves with
oscillations in chemical dynamical systems, and to derive
improved lower bounds for $H_a^c(n)$ and $H_{a,r}^c(n)$,
which are summarize and highlight in red in Table~\ref{tab:2}.
Applying this framework to 
general chemical dynamical systems,
we identify special classes of polynomials, which we call
\emph{semi-chemical} and \emph{chemical} polynomials, 
and investigate the associated Harnack M-curves. 
We use semi-chemical polynomials to construct
a family of $n$-degree chemical dynamical systems with 
$M(n) - (1 + (-1)^n)/2$ algebraic limit cycles for every $n \ge 5$.
Moreover, we construct cubic and quartic 
chemical dynamical systems with $2$ and 
$4$ algebraic limit cycles, respectively, 
which we illustrate in Figure~\ref{fig:1}(a)--(b).
Notably, the cubic example shows that 
chemical dynamical systems can attain 
the best currently known general lower bound for
cubic systems through non-generic algebraic curves.
We then consider chemical dynamical systems admitting 
reversible reaction-network realizations.
In this setting, we identify so-called \emph{buffered} polynomials,
which we use to construct a family of 
reversible $n$-degree chemical dynamical systems with 
$M(n-3)$ algebraic limit cycles for every $n \ge 5$. This establishes a quadratic
lower bound for $H_{a,r}^c(n)$, substantially improving
upon the linear bound summarized in~Table~\ref{tab:1}.
We also construct a family of reversible cubic chemical systems 
with an algebraic limit cycle, and a reversible quartic 
system with $2$ algebraic limit cycles, 
which we illustrate in Figure~\ref{fig:1}(c)--(d).

The remainder of the paper is organized as follows. 
Section~\ref{sec:2} introduces terminology 
and preliminary results concerning chemical dynamical systems and algebraic cycles. 
Section~\ref{sec:3} develops general constructions of planar dynamical systems 
with prescribed algebraic cycles, including perturbation results. 
Sections~\ref{sec:4} and~\ref{sec:5} apply these constructions to 
general and reversible chemical dynamical systems, respectively. 
Finally, Section~\ref{sec:6} discusses the results and presents open problems.

\section{Preliminaries}
\label{sec:2}

\textbf{Notation}. We let 
$\mathbb{Z}_{\ge} = \{0,1,2,\ldots\}$ and
$\mathbb{Z}_{>} = \{1,2,\ldots\}$, and 
denote the non-negative and positive
quadrants of the real plane $\mathbb{R}^2$ by 
$\mathbb{R}_{\ge}^2$ and $\mathbb{R}_{>}^2$, 
respectively. The closure of a set $\Omega \subset \mathbb{R}^N$
is denoted by $\overline{\Omega}$. 

Let $h : \mathbb{R}^2 \to \mathbb{R}$ 
be a polynomial of degree $n \ge 1$, 
which we write as
\begin{align}
h(x,y) & = \sum_{i + j = 0}^n \kappa_{i,j} x^i y^j
= \kappa_{0,0}
+ \sum_{i = 1}^n \kappa_{i,0} x^i
+ \sum_{j = 1}^n \kappa_{0,j} y^j
+ \sum_{\substack{i + j = 2, \\ i,j \ge 1}}^n \kappa_{i,j} x^i y^j.
\label{eq:polynomial_h}
\end{align}

\begin{definition}
$($\textbf{Semi-chemical, chemical and non-chemical polynomials}$)$ 
\label{def:semi_chemical}
Consider an $n$-degree polynomial $h(x,y)$
written as~$(\ref{eq:polynomial_h})$.
If $\kappa_{0,j} \ge 0 \text{ for all } j = 0,1,\ldots,n$, 
then $h(x,y)$ is \emph{chemical in $x$}.
If $\kappa_{i,0} \ge 0 \text{ for all } i = 0,1,\ldots,n$,
then $h(x,y)$ is \emph{chemical in $y$}.
If $h(x,y)$ is chemical in $x$ or $y$, 
then $h$ is \emph{semi-chemical}.
If $h(x,y)$ is chemical in both $x$ and $y$, 
then $h(x,y)$ is \emph{chemical}.
If $h(x,y)$ is neither chemical in $x$ nor in $y$, 
then $h(x,y)$ is \emph{non-chemical}.
\end{definition}

\noindent
\textbf{Remark}. If $h(x,y)$ is chemical in $x$, 
then $p(x,y) = h(y,x)$ is chemical in $y$.

Consider a system of two
autonomous first-order ordinary differential equations 
\begin{align}
\frac{\mathrm{d} x}{\mathrm{d} t} 
& = f(x,y) = \sum_{i + j = 0}^n \alpha_{i,j} x^i y^j
= \sum_{j = 0}^n \alpha_{0,j} y^j
+ \sum_{\substack{i + j = 1, \\ i \ge 1}}^n \alpha_{i,j} x^i y^j, \nonumber \\*
\frac{\mathrm{d} y}{\mathrm{d} t} 
& = g(x,y) = \sum_{i + j = 0}^n \beta_{i,j} x^i y^j
= \sum_{i = 0}^n \beta_{i,0} x^i
+ \sum_{\substack{i + j = 1, \\ j \ge 1}}^n \beta_{i,j} x^i y^j,
\label{eq:DS}
\end{align}
where all the variables are real and dimensionless, 
and $\alpha_{i,j} \ne 0$ or $\beta_{i,j} \ne 0$ 
for some $i + j = n$.
We say that~(\ref{eq:DS}) is a polynomial 
\emph{dynamical system} (DS) of degree $n$~\cite{Perko}. 
The set of all polynomial DSs can be split into 
two disjoint subsets~\cite{QCM}, which we now define.

\begin{definition} 
$($\textbf{Chemical and non-chemical dynamical systems}$)$ 
\label{def:CDS}
Assume that $f(x,y)$ is chemical in $x$,
and that $g(x,y)$ is chemical in $y$. 
Then~$(\ref{eq:DS})$ is called a 
\emph{chemical} dynamical system $($\emph{CDS}$)$.
If~$(\ref{eq:DS})$ is not a \emph{CDS}, 
then it is said to be \emph{non-chemical}.
\end{definition}

Assume that~(\ref{eq:DS}) is a CDS. Then one can associate
with it a set of chemical reactions governed by mass-action kinetics, 
jointly forming a \emph{chemical reaction network} (CRN),
through which the underlying 
chemical species interact~\cite{Feinbergbook,Janos}. 
In particular, let $X$ and $Y$ denote the species with
concentrations $x$ and $y$, respectively.
Then a monomial $x^i y^j$ appearing in~(\ref{eq:DS})
can be associated with a reaction of the form
$i X + j Y \xrightarrow[]{} i' X + j' Y$
for some $i', j' \in \mathbb{Z}_{\ge}$. 
To every CDS one can associate infinitely 
many CRNs~\cite{George_Pantea,NAA}.

\begin{example}
The simple harmonic oscillator
\begin{align}
\frac{\mathrm{d} x}{\mathrm{d} t} 
& = y, \nonumber \\
\frac{\mathrm{d} y}{\mathrm{d} t} 
& = -x, \label{eq:HO}
\end{align}
is a non-chemical \emph{DS} because $g(x,0) = - x$ has negative coefficient,
and is hence not chemical in $y$. 

Consider the system~\emph{$\cite{Escher}$}
\begin{align}
\frac{\mathrm{d} x}{\mathrm{d} t} 
& = \frac{3}{2} + 2 x^2 - x y, \nonumber \\
\frac{\mathrm{d} y}{\mathrm{d} t} 
& = \frac{17}{4} - y + \frac{5}{2} x^2 - x y. \label{eq:Escher}
\end{align}
This is a \emph{CDS}, because $f(0,y) = \frac{3}{2}$
and $g(x,0) = \frac{17}{4} + \frac{5}{2} x^2$
have non-negative coefficients, and are hence chemical
in $x$ and $y$, respectively. 
A \emph{CRN} associated with~$(\ref{eq:Escher})$ is given by
\begin{align}
\varnothing & \xrightarrow[]{\frac{3}{2}} X, \; \; \; 
2 X \xrightarrow[]{2} 3 X,  \; \; \; 
X + Y \xrightarrow[]{1} Y, \nonumber \\
\varnothing & \xrightarrow[]{\frac{17}{4}} Y, \; \; \; 
Y \xrightarrow[]{1} \varnothing, \; \; \; 
2 X \xrightarrow[]{\frac{5}{2}} 2 X + Y, \; \; \; 
X + Y \xrightarrow[]{1} X.
\label{eq:Escher_CRN}
\end{align}
Here, the \emph{rate coefficients}, indicated 
above the reaction arrows, determine the rates 
of the corresponding reactions under mass-action kinetics.
\end{example}

\subsection{Algebraic cycles and their stability} 
In what follows, we define \emph{period} $\tau$
of a non-constant continuous periodic function 
$\mathbf{x} : \mathbb{R} \to \mathbb{R}^2$ 
as the smallest number $\tau > 0$ such that 
$\mathbf{x}(t + \tau) = \mathbf{x}(t)$ 
for every $t \in \mathbb{R}$. 

\begin{definition} $($\textbf{Cycle}$)$ 
\label{def:cycle}
Let $\mathbf{x}_{\Gamma} : \mathbb{R} \to \mathbb{R}^2$,
$\mathbf{x}_{\Gamma}(t) = (x(t),y(t))$,
be a non-constant periodic solution of~$(\ref{eq:DS})$ 
with period $\tau$. 
Then the set $\Gamma = \big\{\mathbf{x}_{\Gamma}(t)
\,\big|\, t \in [0,\tau) \big\}$ is called a \emph{cycle} of~$(\ref{eq:DS})$.
If there exists an open set
$\mathbb{U} \supset \Gamma$
such that $\Gamma$ is the only
cycle contained in $\mathbb{U}$, 
then $\Gamma$ is called a \emph{limit cycle}. 
\end{definition}

\noindent
Given a non-constant polynomial $h : \mathbb{R}^2 \to \mathbb{R}$,
we denote its zero-level set by 
$\{h = 0\} = \big\{(x,y) \in \mathbb{R}^2 \,\big|\,
h(x,y) = 0\big\}$. This set is called an algebraic curve.

\begin{definition} $($\textbf{Algebraic cycle}$)$ 
\label{def:a_cycle}
Let $\Gamma$ be a cycle of~$(\ref{eq:DS})$.
If there exists a non-constant polynomial $h : \mathbb{R}^2 \to \mathbb{R}$
such that $\Gamma \subset \{h = 0\}$, 
then $\Gamma$ is called an 
\emph{algebraic cycle} of~$(\ref{eq:DS})$.
\end{definition}

Let $\nabla h = (\frac{\partial h}{\partial x},
\frac{\partial h}{\partial y})$
denote the gradient of $h$.
We now define a notion of invariance.

\begin{definition}
\label{def:cofactor}
$($\textbf{Invariant algebraic curve and cofactor}$)$ 
Let $h : \mathbb{R}^2 \to \mathbb{R}$ be a non-constant polynomial.
The algebraic curve $\{h = 0\}$ is said to be \emph{invariant}
with respect to \emph{DS}~$(\ref{eq:DS})$ if there exists 
a polynomial $k : \mathbb{R}^2 \to \mathbb{R}$ such that 
\begin{align}
(\mathbf{f} \cdot \nabla h)(x,y) 
& = f(x,y) \frac{\partial h(x,y)}{\partial x}
+ g(x,y) \frac{\partial h(x,y)}{\partial y}
 = k(x,y) h(x,y)
\; \; \; 
\text{for every } (x,y) \in \mathbb{R}^2.
\label{eq:cofactor_def}
\end{align}
The polynomial $k$ is called the \emph{cofactor} of $h$ 
with respect to \emph{DS}~$(\ref{eq:DS})$.
\end{definition}

\noindent 
\textbf{Remark}. It follows directly from~(\ref{eq:cofactor_def})
that, if DS~(\ref{eq:DS}) is of degree $n$, then 
the cofactor $k$, if non-zero, is of degree at most $n - 1$. 
Furthermore, the cofactor is unique. 
Indeed, if~(\ref{eq:cofactor_def}) holds with two 
cofactors $k_1$ and $k_2$, then $(k_1 - k_2) h \equiv 0$.
Since $h \not \equiv 0$, it follows that 
$k_1 = k_2$ on $\mathbb{R}^2$, proving uniqueness.

\noindent
\textbf{Remark}. The notion of invariance
introduced in Definition~\ref{def:cofactor} is stronger
than the usual dynamical-systems notion of invariance.
Indeed, if $\{h = 0\}$ is invariant 
in the sense of Definition~\ref{def:cofactor}, 
then $\frac{\mathrm{d} h}{\mathrm{d} t} = k(x,y) h$
along every solution of DS~(\ref{eq:DS}).
Consequently, any trajectory starting in 
$\{h = 0\}$ remains in this set, so 
$\{h = 0\}$ is also invariant in the usual sense.
However, the converse is not true in general 
for a given polynomial $h$, since the existence 
of a polynomial cofactor is not guaranteed.

\vskip 2mm

\noindent \textbf{Hyperbolicity and stability}.
Let $\mathbf{x}(t;\mathbf{x}_0) \in \mathbb{R}^2$
be the solution of~(\ref{eq:DS}) satisfying
$\mathbf{x}(0;\mathbf{x}_0) = \mathbf{x}_0$.
Given a $\tau$-periodic cycle $\Gamma$, we define
the distance from a point $\mathbf{x}(t)$ to $\Gamma$ by
\begin{align}
\textrm{dist}(\mathbf{x}(t),\Gamma) & = \textrm{min}_{s \in [0,\tau]} 
\sqrt{(x(t) - x_{\Gamma}(s))^2 + (y(t) - y_{\Gamma}(s))^2}.
\end{align}

\begin{definition} $($\textbf{Hyperbolicity and stability}$)$ 
\label{def:stability}
A limit cycle $\Gamma$ of~$(\ref{eq:DS})$
is said to be \emph{hyperbolic} if 
there exists an open set 
$\mathbb{U} \supset \Gamma$ and constants
$a, \alpha > 0$ such that one of the following
two conditions holds.
\begin{enumerate}
\item[] \textbf{Stable case}. For every 
$\mathbf{x}_0 \in \mathbb{U}$ and $t \ge 0$,
\begin{align}
\textrm{\emph{dist}}(\mathbf{x}(t;\mathbf{x}_0),\Gamma)
\le a \textrm{\emph{dist}}(\mathbf{x}_0,\Gamma) \exp(-\alpha t). \label{eq:stable}
\end{align}
\item[] \textbf{Unstable case}. For every 
$\mathbf{x}_0 \in \mathbb{U}$ and $t \le 0$,
\begin{align}
\textrm{\emph{dist}}(\mathbf{x}(t;\mathbf{x}_0),\Gamma)
\le a \textrm{\emph{dist}}(\mathbf{x}_0,\Gamma) \exp(\alpha t).
\label{eq:unstable}
\end{align}
\end{enumerate}
A hyperbolic limit cycle $\Gamma$ is said to be
\emph{stable} if~$(\ref{eq:stable})$ holds
and \emph{unstable} if~$(\ref{eq:unstable})$ holds.
\end{definition} 

A connected component 
$\Gamma \subset \{h = 0\}$ is said to be 
\emph{regular} (non-singular)  
if $\nabla h(x,y) \ne \mathbf{0}$
for every $(x,y) \in \Gamma$.
The hyperbolicity and stability of a regular
cycle contained in an invariant 
algebraic curve can be determined from 
a period integral of its cofactor. 

\begin{theorem}
\label{theorem:cofactor}
Assume that \emph{DS}~$(\ref{eq:DS})$ admits an 
invariant algebraic curve $\{h = 0\}$ with polynomial cofactor $k$.
Let $\Gamma \subset \{h = 0\}$ be a regular algebraic 
cycle with period $\tau$, and define
\begin{align}
K_{\Gamma} & = \int_0^{\tau} 
k(\mathbf{x}_{\Gamma}(t)) \mathrm{d} t.
\label{eq:K_cofactor}
\end{align}
Then $\Gamma$ is a hyperbolic limit cycle if and only if
$K_{\Gamma} \ne 0$. $\Gamma$ is stable if $K_{\Gamma} < 0$
and unstable if $K_{\Gamma} > 0$.
\end{theorem}

\begin{proof}
Since $\nabla h \ne \mathbf{0}$ on $\Gamma$, 
by continuity, there exists an open set $\mathbb{U} \supset \Gamma$
on which $\nabla h \ne \mathbf{0}$. 
Since $k$ is a polynomial cofactor of $h$ on $\mathbb{R}^2$,
the assumptions of~\cite{Giacomini}[Theorem 2] hold, giving
$\int_0^{\tau} k(\mathbf{x}_{\Gamma}(t)) \mathrm{d} t 
= \int_0^{\tau} \left(\partial f/\partial x
+ \partial g/\partial y \right)(\mathbf{x}_{\Gamma}(t)) \mathrm{d} t$,
from which the conclusions follow.
\end{proof}

\section{Dynamical systems with prescribed algebraic cycles}
\label{sec:3}
Let us consider the class of polynomial DSs~(\ref{eq:DS}) taking the following form:
\begin{align}
\frac{\mathrm{d} x}{\mathrm{d} t} 
& = f(x,y) = 
f_1(x,y) h(x,y) - c(x,y) \frac{\partial h(x,y)}{\partial y}, 
\nonumber \\
\frac{\mathrm{d} y}{\mathrm{d} t} 
& = g(x,y) = 
g_1(x,y) h(x,y) + c(x,y) \frac{\partial h(x,y)}{\partial x},
\label{eq:DS_1}
\end{align} 
where $f_1, g_1, c, h : \mathbb{R}^2 \to \mathbb{R}$ are polynomials.
The class~(\ref{eq:DS_1}) encompasses a broad range of polynomial
DSs~(\ref{eq:DS}) admitting $\{h = 0\}$ as an invariant algebraic curve. 
In particular, a sufficient condition for DS~(\ref{eq:DS}) 
possessing $\{h = 0\}$ as an invariant 
algebraic curve, in the sense of Definition~\ref{def:cofactor}, 
to admit such a representation 
is that $\{h = 0\}$ is regular in $\mathbb{C}^2$~\cite{Inverse}[Theorem 3.8(a)], 
where $\mathbb{C}$ denotes the set of complex numbers. 
All DSs constructed in this paper are of the form~(\ref{eq:DS_1}).

\begin{example}
Ellipse $\{h = 0\}$ is an invariant algebraic curve 
of \emph{CDS}~$(\ref{eq:Escher})$~\emph{\cite{Escher}}, where
\begin{align}
h(x, y) & = 19 + 20 x - 16 y + 10 x^2 + 4 y^2 - 12 x y. 
\label{eq:Escher_h}
\end{align} 
The unique solution of $\nabla h = \mathbf{0}$ in $\mathbb{C}^2$ 
is $(2,5)$. Since $h(2,5) = - 1 \ne 0$, the algebraic curve $\{h = 0\}$ is 
regular in $\mathbb{C}^2$. Therefore, $(\ref{eq:Escher})$
can be written in the form~$(\ref{eq:DS_1})$. Explicit computation gives
\begin{align}
f_1(x,y) & = \frac{1}{2}, \; \; \; 
g_1(x,y) = \frac{3}{4}, \; \; \; 
c(x,y) = \frac{1}{4} (-2 - x + y).
\label{eq:Escher_functions}
\end{align}
\end{example}

In~\cite{Winkel}, the focus was on DSs~(\ref{eq:DS_1}) in the special case where 
$f_1 = - \frac{\partial h}{\partial x}$, $g_1 = - \frac{\partial h}{\partial y}$ and $c = 1$,
while~\cite{Christopher} considered the case where $f_1$ and $g_1$ are constants and $c$ is linear. 
We now generalize these results by allowing $f_1,g_1$ and $c$ 
to be arbitrary polynomials. To this end, for every $\tau$-periodic cycle 
$\Gamma$ of~(\ref{eq:DS_1}), define
\begin{align}
K_{\Gamma} = 
\int_0^{\tau} 
\left(f_1(\mathbf{x}_{\Gamma}(t)) 
\frac{\partial h(\mathbf{x}_{\Gamma}(t))}{\partial x} 
+ g_1(\mathbf{x}_{\Gamma}(t)) 
\frac{\partial h(\mathbf{x}_{\Gamma}(t))}{\partial y} \right) \mathrm{d} t.
\label{eq:K}
\end{align}

\begin{theorem}
\label{theorem:AL}
Consider polynomial \emph{DS}~$(\ref{eq:DS_1})$.
The set $\{h = 0\}$ is an invariant algebraic curve of~$(\ref{eq:DS_1})$.
Let $\Gamma \subset \{h = 0\}$ 
be a regular simple closed curve,
and assume that $c(x,y) \ne 0$ 
for every $(x,y) \in \Gamma$. 
Then $\Gamma$ is an algebraic cycle of 
\emph{DS}~$(\ref{eq:DS_1})$. Furthermore, 
denoting its period by $\tau$,
$\Gamma$ is a hyperbolic limit cycle 
of~$(\ref{eq:DS_1})$ if and only if 
$K_{\Gamma} \ne 0$, where $K_{\Gamma}$ is given by~$(\ref{eq:K})$.
In this case, $\Gamma$ is stable if $K_{\Gamma} < 0$
and unstable if $K_{\Gamma} > 0$.
\end{theorem}

\begin{proof}
Differentiating $h = h(x(t),y(t))$
with respect to $t$ gives
\begin{align}
\frac{\mathrm{d} h}{\mathrm{d} t}
& = \frac{\partial h}{\partial x} \left(f_1 h
- c \frac{\partial h}{\partial y} \right)
+ \frac{\partial h}{\partial y} 
\left(g_1 h + c \frac{\partial h}{\partial x} \right) = k h,
\label{eq:h_invariance}
\end{align} 
where $k(x,y)$ is the polynomial cofactor of 
$h$ with respect to DS~(\ref{eq:DS_1}) and reads as
\begin{align}
k(x,y) & = f_1(x,y) \frac{\partial h(x,y)}{\partial x}
+  g_1(x,y) \frac{\partial h(x,y)}{\partial y}.
\label{eq:cofactor_1}
\end{align}
Hence, by Definition~(\ref{def:cofactor}), 
$\{h = 0\}$ is an invariant algebraic curve of~$(\ref{eq:DS_1})$.
In particular, $\Gamma$ is an invariant simple closed curve.
Every equilibrium $(x^*,y^*) \in \Gamma$ satisfies
\begin{align}
0 & = c(x^*,y^*) \frac{\partial h(x^*,y^*)}{\partial y}, \; \; 
0 = c(x^*,y^*) \frac{\partial h(x^*,y^*)}{\partial x}.
\label{eq:equilibria_AL}
\end{align}
Since $\nabla h(x^*,y^*) \ne \mathbf{0}$ and $c(x^*,y^*) \ne 0$, 
(\ref{eq:equilibria_AL}) has no solutions.
Hence $\Gamma$ is an invariant simple closed
curve along which the vector field does not vanish, 
and is therefore an algebraic cycle.
Substituting the cofactor~(\ref{eq:cofactor_1})
into~(\ref{eq:K_cofactor}), the hyperbolicity 
and stability results follow 
from Theorem~\ref{theorem:cofactor}.
\end{proof}

\noindent
\textbf{Remark}. Under the hypotheses of Theorem~\ref{theorem:AL}, if $f_1 = - \frac{\partial h}{\partial x}$
and $g_1 = - \frac{\partial h}{\partial y}$,
then $K_{\Gamma} = - \int_0^{\tau} 
\big( \big(\frac{\partial h}{\partial x} \big)^2
+ \big(\frac{\partial h}{\partial y} \big)^2 \big)
\mathrm{d} t < 0$. Hence, in this case, 
every regular simple closed curve $\Gamma
\subset \{h = 0\}$ is a stable limit cycle of~(\ref{eq:DS_1}); 
see also~\cite{Winkel}. 

\vskip 2mm

\noindent
Theorem~\ref{theorem:AL} assumes that $c \ne 0$ on $\Gamma$.
If, additionally, $c \ne 0$ in the interior of $\Gamma$,
then one can express $K_{\Gamma}$ from~(\ref{eq:K}),
whose integrand depends on $h$, with an integral
whose integrand is independent of $h$. 
In what follows, the sign function 
$\textrm{sign} : \mathbb{R} \to \{-1,0,1\}$
is defined as
$\textrm{sign}(x) = -1$ if $x < 0$, 
$\textrm{sign}(x) = 0$ if $x = 0$, and 
$\textrm{sign}(x) = 1$ if $x > 0$.

\begin{lemma}
\label{lemma:AL_Green}
Let $\Gamma$ be a cycle of \emph{DS}~$(\ref{eq:DS_1})$. 
Assume that $c(x,y) \ne 0$ for every $(x,y) \in
\overline{\Omega(\Gamma)}$, 
where $\Omega(\Gamma)$ is the interior of $\Gamma$. 
Then $K_{\Gamma}$, given by~$(\ref{eq:K})$, satisfies 
\begin{align}
K_{\Gamma} & = s(\mathbf{x}_{\Gamma}) 
\iint_{\Omega(\Gamma)} \frac{F(x,y)}{c(x,y)^2}
\mathrm{d} x \mathrm{d} y,
\label{eq:K_F}
\end{align}
where 
\begin{align}
F & = c \left(\frac{\partial f_1}{\partial x} 
+ \frac{\partial g_1}{\partial y}\right)
- \left(f_1 \frac{\partial c}{\partial x} 
+ g_1 \frac{\partial c}{\partial y} \right),
\label{eq:F}
\end{align}
and $s(\mathbf{x}_{\Gamma}) = + 1$ (respectively, $s(\mathbf{x}_{\Gamma}) = - 1$)
if $\mathbf{x}_{\Gamma}(t)$ is traversed 
counter-clockwise (respectively, clockwise) as $t$ increases.
If, furthermore, $\Gamma \subset \{h = 0\}$, then
\begin{align}
s(\mathbf{x}_{\Gamma}) & = 
\textrm{\emph{sign}}(c(\mathbf{p}))
\textrm{\emph{sign}}(h(\mathbf{q}))  \label{eq:s}
\end{align}
for every point $\mathbf{p} \in \Gamma$,
and every point $\mathbf{q}$
sufficiently close to $\Gamma$ and in 
the exterior of $\Gamma$.
\end{lemma}

\begin{proof}
Since $c \ne 0$ on $\Gamma$, multiplying the first 
and second equations in~(\ref{eq:DS_1}) by $g_1/c$ and $f_1/c$,
respectively, and subtracting the former from the latter, 
one obtains along $\Gamma$
\begin{align}
\frac{f_1}{c} \frac{\mathrm{d} y}{\mathrm{d} t} 
- \frac{g_1}{c} \frac{\mathrm{d} x}{\mathrm{d} t}
& = f_1 \frac{\partial h}{\partial x} + g_1 \frac{\partial h}{\partial y}.
\end{align}
Hence, by~(\ref{eq:K}), 
\begin{align}
K_{\Gamma} 
& =  \int_0^{\tau}
\left(\frac{f_1}{c}  \frac{\mathrm{d} y}{\mathrm{d} t}
- \frac{g_1}{c}  \frac{\mathrm{d} x}{\mathrm{d} t} \right) \mathrm{d} t
= \int_0^{\tau} 
\left(\frac{\mathbf{f}_1(\mathbf{x}_{\Gamma}(t))}{
c(\mathbf{x}_{\Gamma}(t))}
\right) 
\cdot \mathbf{N}(\mathbf{x}_{\Gamma}(t)) \mathrm{d} t,
\end{align}
where $\mathbf{f}_1 = (f_1,g_1)$, and 
$\mathbf{N} = (\frac{\mathrm{d} y}{\mathrm{d} t},-\frac{\mathrm{d} x}{\mathrm{d} t})$
evaluated at $\mathbf{x}_{\Gamma}(t)$
is a right-pointing normal vector to
$\mathbf{x}_{\Gamma}(t)$. 
By the assumption that $c(x,y) \ne 0$ 
for every $(x,y) \in \overline{\Omega(\Gamma)}$,
and since $\overline{\Omega(\Gamma)}$ is compact and $c$ is continuous, 
it follows that $c(x,y) \ne 0$ 
on some open set $\mathbb{U} \supset 
\overline{\Omega(\Gamma)}$;
hence, $\mathbf{f}_1/c$ is continuously differentiable on $\mathbb{U}$. 
Therefore, Green's theorem yields
\begin{align}
K_{\Gamma} 
& = \int_0^{\tau} 
\left(\frac{\mathbf{f}_1(\mathbf{x}_{\Gamma}(t))}{
c(\mathbf{x}_{\Gamma}(t))} \right)
\cdot \mathbf{N}(\mathbf{x}_{\Gamma}(t)) \mathrm{d} t
= s(\mathbf{x}_{\Gamma}) \iint_{\Omega(\Gamma)} 
\nabla \cdot \left(\frac{\mathbf{f}_1(x,y)}{c(x,y)} \right) 
\mathrm{d} x \mathrm{d} y. \label{eq:Green}
\end{align}
Expanding the right-hand side in~(\ref{eq:Green}), 
one obtains~(\ref{eq:K_F}).

Assume now that $\Gamma \subset \{h = 0\}$. Then, along $\Gamma$, 
\begin{align}
\mathbf{N}(\mathbf{x}_{\Gamma}(t)) & = \left(\frac{\mathrm{d} y_{\Gamma}}{\mathrm{d} t}, 
- \frac{\mathrm{d} x_{\Gamma}}{\mathrm{d} t} \right) = 
c(\mathbf{x}_{\Gamma}(t)) \nabla h(\mathbf{x}_{\Gamma}(t)).
\label{eq:normal}
\end{align}
Assume first that $\mathbf{x}_{\Gamma}(t)$ is traversed counter-clockwise
as $t$ increases. Then $s(\mathbf{x}_{\Gamma}) = +1$
and $\mathbf{N}(\mathbf{x}_{\Gamma})$ points outwards. 
If $c(\mathbf{x}_{\Gamma}) > 0$, then $\nabla h(\mathbf{x}_{\Gamma})$ 
points outwards by~(\ref{eq:normal}), 
and hence $h > 0$ immediately outside $\Gamma$. 
If $c(\mathbf{x}_{\Gamma}) < 0$, then $\nabla h(\mathbf{x}_{\Gamma})$ 
points inwards, so $h < 0$ immediately outside $\Gamma$. Hence, for
any $\mathbf{p} \in \Gamma$ and any $\mathbf{q} \in \mathbb{R}^2$
in the immediately exterior of $\Gamma$, 
$c(\mathbf{p})$ and $h(\mathbf{q})$ have the same sign,
and $\textrm{sign}(c(\mathbf{p})) \textrm{sign}(h(\mathbf{q}))  = + 1 = s$.
If $\mathbf{x}_{\Gamma}(t)$ is traversed clockwise, 
then $\mathbf{N}(\mathbf{x}_{\Gamma})$ points inwards. 
The same reasoning then shows that $\textrm{sign}(c(\mathbf{p})) 
\textrm{sign}(h(\mathbf{q}))  = - 1 = s$, thus proving~(\ref{eq:s}).
\end{proof}

We now generalize~\cite{Christopher}[Theorem 1].

\begin{theorem}
\label{theorem:AL_strong}
Consider polynomial \emph{DS}~$(\ref{eq:DS_1})$.
Let $\mathbb{S} \subseteq \mathbb{R}^2$ be 
an open simply connected set such that, for every $(x,y) \in \mathbb{S}$, 
$c(x,y) \ne 0$ and $F(x,y) \ne 0$, where $F$ is given by~$(\ref{eq:F})$.
Assume that $\{h = 0\}$ contains 
$N$ regular simple closed curves 
$\Gamma_1, \Gamma_2, \ldots, \Gamma_N \subset \mathbb{S}$. 
Then $\Gamma_1, \Gamma_2, \ldots, \Gamma_N$
are hyperbolic algebraic limit cycles of \emph{DS}~$(\ref{eq:DS_1})$. 
Moreover, $\Gamma_i$ is stable (respectively, unstable)
if $K_{\Gamma_i} < 0$ (respectively, $K_{\Gamma_i} > 0$),
where $K_{\Gamma_i}$ is given by~$(\ref{eq:K_F})$
with $s(\mathbf{x}_{\Gamma_i})$ satisfying~$(\ref{eq:s})$.
Furthermore, if $\{h = 0\}$ contains no 
bounded connected components in $\mathbb{S}$ other than 
$\Gamma_1, \Gamma_2, \ldots, \Gamma_N$, 
then~$(\ref{eq:DS_1})$ has no other cycles in $\mathbb{S}$.
\end{theorem}

\begin{proof}
Since $\mathbb{S}$ is simply connected, 
for every $i = 1,2,\ldots, N$, $\overline{\Omega(\Gamma_i)}
\subset \mathbb{S}$ and hence $c(x,y) \ne 0$ 
for every $(x,y) \in \overline{\Omega(\Gamma_i)}$.
Therefore, by Theorem~\ref{theorem:AL}, 
$\Gamma_1, \Gamma_2, \ldots, \Gamma_N$ are cycles of~(\ref{eq:DS_1})
and Lemma~\ref{lemma:AL_Green} applies.
Since $\mathbb{S}$ is connected, $F$ is continuous 
and $F \ne 0$ on $\mathbb{S}$, it follows
that $F$ has constant sign on $\mathbb{S}$
and $K_{\Gamma_i} \ne 0$.
Therefore, by Lemma~\ref{lemma:AL_Green}, 
$\Gamma_1, \Gamma_2, \ldots, \Gamma_N$ 
are hyperbolic limit cycles of~(\ref{eq:DS_1})
with stated stability.

Assume that $\{h = 0\}$ contains no 
bounded connected components in $\mathbb{S}$ other than 
$\Gamma_1,$ $\Gamma_2,$ $\ldots,$ $\Gamma_N$. 
To prove that~(\ref{eq:DS_1}) 
then has no other cycle in $\mathbb{S}$, 
assume for contradiction that 
there exists a cycle $\tilde{\Gamma} \subset \mathbb{S}$
such that $\tilde{\Gamma} \not\in 
\{\Gamma_1, \Gamma_2, \ldots, \Gamma_N\}$.
Let $\tilde{\tau}$ and $\Omega(\tilde{\Gamma}) 
\subset \mathbb{S}$ be the period and the interior of 
$\tilde{\Gamma}$, respectively. 

Assume for contradiction that $\tilde{\Gamma}$ intersects $\{h = 0\}$.
By invariance of $\{h = 0\}$ and uniqueness of solutions, 
it then follows that $\tilde{\Gamma} \subset \{h = 0\}$.
Hence, by~(\ref{eq:DS_1}), the vector field satisfies 
$\|\mathbf{f}(x,y)\|^2 = c(x,y)^2 \|\nabla h(x,y)\|^2$ on $\tilde{\Gamma}$. 
Since the vector field does not vanish on $\tilde{\Gamma}$, it follows 
that $\nabla h \ne \mathbf{0}$ on $\tilde{\Gamma}$. 
Consequently, $\tilde{\Gamma}$ is a regular bounded
connected component of $\{h = 0\}$, contradicting our assumption.
Therefore, $h(x,y) \ne 0$ for every $(x,y) \in \tilde{\Gamma}$. 

Since $h \ne 0$ on $\tilde{\Gamma}$, dividing both sides of~(\ref{eq:h_invariance})
by $h$ gives, along $\tilde{\Gamma}$, 
\begin{align}
k & = \frac{1}{h} \frac{\mathrm{d} h}{\mathrm{d} t} = \frac{\mathrm{d}}{\mathrm{d} t} \log |h|.
\end{align}
By~(\ref{eq:K}) and~(\ref{eq:cofactor_1}),
\begin{align}
K_{\tilde{\Gamma}} & = \int_0^{\tilde{\tau}} 
k(\mathbf{x}_{\tilde{\Gamma}}(t)) \mathrm{d} t
 = \int_0^{\tilde{\tau}} \frac{\mathrm{d}}{\mathrm{d} t} \log|h(\mathbf{x}_{\tilde{\Gamma}}(t))| \mathrm{d} t = 0, \label{eq:Green_LHS}
\end{align}
where the last equality follows from the periodicity of $\mathbf{x}_{\tilde{\Gamma}}$. 
On the other hand, since $\mathbb{S}$ is simply connected, 
$\overline{\Omega(\tilde{\Gamma})} \subset \mathbb{S}$.
Since $c \ne 0$ on $\mathbb{S}$, and~(\ref{eq:K_F}) holds also 
for cycles not contained in $\{h = 0\}$, Lemma~\ref{lemma:AL_Green} yields
\begin{align}
K_{\tilde{\Gamma}} & = s(\mathbf{x}_{\tilde{\Gamma}}) \iint_{\Omega(\tilde{\Gamma})} 
\frac{F(x,y)}{c(x,y)^2} \mathrm{d} x \mathrm{d} y. \label{eq:Green_RHS}
\end{align}
Since $F$ has constant non-zero sign on $\mathbb{S}$,
it follows that $K_{\tilde{\Gamma}} \ne 0$, contradicting~(\ref{eq:Green_LHS}).
Therefore, there are no cycles in $\mathbb{S}$ other than 
$\Gamma_1, \Gamma_2, \ldots, \Gamma_N$. 
\end{proof}

\noindent 
\textbf{Remark}. To exclude the existence
of cycles other than $\Gamma_1, \Gamma_2, \ldots, \Gamma_N$
globally in $\mathbb{R}^2$, a quantity of interest 
is the derivative of $c$ along the solutions on $\{c = 0\}$,
which is given by 
\begin{align}
\frac{\mathrm{d} c}{\mathrm{d} t} \big|_{\{c = 0\}}
& = - F h \big|_{\{c = 0\}}. \label{eq:c_derivative}
\end{align}
If $F h$ has a constant non-zero sign 
on $\{c = 0\}$, then no cycle can intersect $\{c = 0\}$.
Furthermore, if every connected component of $\{c = 0\}$
is unbounded, then $\{c = 0\}$ cannot intersect the
interior of any cycle. 

Let us apply Theorem~\ref{theorem:AL_strong} to show that
CDS~(\ref{eq:Escher}) has a stable hyperbolic 
algebraic limit cycle in $\mathbb{R}_{>}^2$, 
as established by construction in~\cite{Escher}.
We additionally prove that this cycle is unique in $\mathbb{R}^2$. 

\begin{lemma}
\label{lemma:Escher}
\emph{CDS}~$(\ref{eq:Escher})$ has a stable hyperbolic
algebraic limit cycle in $\mathbb{R}_{>}^2$, given as the ellipse
$\{h = 0\}$, where $h$ is given by~$(\ref{eq:Escher_h})$.
Furthermore, $(\ref{eq:Escher})$ has no other cycles in $\mathbb{R}^2$.
\end{lemma}

\begin{proof}
CDS~$(\ref{eq:Escher})$ is of the form~$(\ref{eq:DS_1})$ with
$f_1, g_1, c$ given in~(\ref{eq:Escher_functions}). 
The zero-level set $\{h = 0\}$ is a regular ellipse $\Gamma$ contained
in $\mathbb{R}_{>}^2$. Furthermore, the function $F$, defined in~(\ref{eq:F}), 
satisfies $F(x,y) = -1/16 < 0$ on $\mathbb{R}^2$. 
Consider the open simply connected sets
$\mathbb{S}_{+} = \{(x,y) \in \mathbb{R}^2 \, | \, y > x + 2\}$ and 
$\mathbb{S}_{-} = \{(x,y) \in \mathbb{R}^2 \, | \, y < x + 2\}$.
Since $c > 0$ on $\mathbb{S}_{+}$ and $c < 0$ on $\mathbb{S}_{-}$, 
the assumptions about $F$ and $c$ from Theorem~\ref{theorem:AL_strong}
hold on both sets. Furthermore, $h(x, x + 2) = 2 (x - 1)^2 + 1 > 0$,
so the line $y = x + 2$ does not intersect $\Gamma$. 
Since the center $(2,5)$ of $\Gamma$ lies in $y > x + 2$, 
it follows that $\Gamma\subset \mathbb{S}_{+}$. Since $h > 0$ outside $\Gamma$, 
and $F < 0$ and $c > 0$ in the region enclosed by $\Gamma$,  
$K_{\Gamma} < 0$ by~(\ref{eq:K_F}). 
Therefore, by Theorem~\ref{theorem:AL_strong}, $\Gamma$ is a 
hyperbolic algebraic limit cycle which is stable, 
and there are no other cycles in $\mathbb{S}_{+}$ 
and no cycles in $\mathbb{S}_{-}$.
Furthermore, by~(\ref{eq:c_derivative}),
\begin{align}
\frac{\mathrm{d}}{\mathrm{d} t} c \big|_{y = 2 + x}
& = \frac{1}{16} h(x, x + 2) = \frac{1}{16} 
\left(2 (x - 1)^2 + 1 \right) > 0.
\end{align}
Therefore, no cycle intersects $\{c = 0\}$. 
Hence, $\Gamma$ is the unique cycle in $\mathbb{R}^2$.
\end{proof}

\noindent
\textbf{Remark}. Lemma~\ref{lemma:Escher} is consistent 
with the known result that a quadratic DS admitting an invariant 
algebraic curve of degree at most two can have a limit cycle
only if the curve is an ellipse (in which case the limit cycle
coincides with the ellipse), a parabola or a straight line.
Moreover, in each of these cases, if a limit cycle exists, 
it is unique and hyperbolic; see, for example, 
\cite{Gasull_Giacomini,Zegeling} 
and the references therein.

\subsection{Perturbations preserving algebraic limit cycles}
Consider the polynomial DS
\begin{align}
\frac{\mathrm{d} x}{\mathrm{d} t} 
& = f_0(x,y), \nonumber \\
\frac{\mathrm{d} y}{\mathrm{d} t} 
& = g_0(x,y),
\label{eq:DS_0}
\end{align}
and assume that $\{h = 0\}$ is an invariant algebraic curve
of this system containing a hyperbolic algebraic limit cycle $\Gamma$.
By hyperbolicity, this limit cycle persists
under every sufficiently small polynomial 
perturbation of~(\ref{eq:DS_0})~\cite{Perko}. 
However, the perturbed limit cycle need not be algebraic.  
In this section, we consider a special class of perturbations
that preserve $\Gamma$ itself as an algebraic limit cycle.
In particular, let us perturb~(\ref{eq:DS_0})
by adding the right-hand side of~(\ref{eq:DS_1})
multiplied by a parameter $\mu \in \mathbb{R}$, with $|\mu|$ small,
\begin{align}
\frac{\mathrm{d} x}{\mathrm{d} t} 
& = f(x,y; \mu) = f_0(x,y) + \mu \left(f_1(x,y) h(x,y) 
- c(x,y) \frac{\partial h(x,y)}{\partial y} \right), 
\nonumber \\
\frac{\mathrm{d} y}{\mathrm{d} t} 
& = g(x,y; \mu) = g_0(x,y) + \mu \left(g_1(x,y) h(x,y) 
+ c(x,y) \frac{\partial h(x,y)}{\partial x} \right).
\label{eq:DS_0_mu}
\end{align} 

\begin{theorem}
\label{theorem:perturbations}
Assume that $\{h = 0\}$ is an invariant algebraic curve
of polynomial \emph{DS}~$(\ref{eq:DS_0})$.
Then $\{h = 0\}$ is also an invariant algebraic curve
of polynomial \emph{DS}~$(\ref{eq:DS_0_mu})$ 
for every $\mu \in \mathbb{R}$.
Furthermore, if $\Gamma \subset \{h = 0\}$ is a 
hyperbolic algebraic limit cycle of~$(\ref{eq:DS_0})$,
then there exists $\mu_0 > 0$ such that, 
for every $\mu \in (-\mu_0,\mu_0)$, 
$\Gamma$ is a hyperbolic algebraic 
limit cycle of \emph{DS}~$(\ref{eq:DS_0_mu})$
with the same stability as in~$(\ref{eq:DS_0})$.
\end{theorem}

\begin{proof}
By Definition~\ref{def:cofactor}, 
there exists a polynomial cofactor $k_0$ of $h$
with respect to DS~(\ref{eq:DS_0}),
\begin{align}
f_0(x,y) \frac{\partial h(x,y)}{\partial x} +  
g_0(x,y) \frac{\partial h(x,y)}{\partial y} 
& = k_0(x,y) h(x,y) 
\; \; \; \text{for every } (x,y) \in \mathbb{R}^2.
\label{eq:h_invariance_0}
\end{align}
Computing the derivative of $h$ along the solutions 
of~(\ref{eq:DS_0_mu}), one obtains
\begin{align}
\frac{\mathrm{d} h}{\mathrm{d} t}
& = f \frac{\partial h}{\partial x}
+ g \frac{\partial h}{\partial y}
= \left(f_0 \frac{\partial h}{\partial x}
+ g_0 \frac{\partial h}{\partial y} \right)
+ \mu \left(f_1 \frac{\partial h}{\partial x}
+  g_1 \frac{\partial h}{\partial y} \right) h.
\end{align} 
Substituting~(\ref{eq:h_invariance_0}), it follows that
\begin{align}
f(x,y;\mu) \frac{\partial h(x,y)}{\partial x} 
+ g(x,y;\mu) \frac{\partial h(x,y)}{\partial y}
& = k_{\mu}(x,y) h(x,y)
\; \; \text{for every } (x,y) \in \mathbb{R}^2,
\end{align} 
where the cofactor of $h$
with respect to DS~(\ref{eq:DS_0_mu}) is given by
\begin{align}
k_{\mu}(x,y) & = k_0(x,y) + \mu \left(f_1(x,y)
\frac{\partial h(x,y)}{\partial x} + g_1(x,y) 
\frac{\partial h(x,y)}{\partial y} \right).
\label{eq:h_invariance_mu}
\end{align}
Therefore $\{h = 0\}$ is an invariant algebraic curve
of~(\ref{eq:DS_0_mu}) for all $\mu \in \mathbb{R}$.
In particular, $\Gamma$ is an invariant simple closed curve.

On $\Gamma$, the perturbed vector field satisfies
\begin{align}
\mathbf{f}(x,y; \mu) & = \mathbf{f}_0(x,y) + \mu c(x,y) 
\left( -\frac{\partial h}{\partial y}, \frac{\partial h}{\partial x} \right).
\end{align}
Since $\Gamma$ is a cycle of~(\ref{eq:DS_0}) and is compact,
$m = \min_{(x,y) \in \Gamma} \|\mathbf{f}_0(x,y)\| > 0$. 
Furthermore, by continuity, 
$M = \max_{(x,y) \in \Gamma} |c(x,y)| \|\nabla h(x,y)\| < \infty$.
Choose $\mu_0 = m/(2 M + 1)$. Then, for every $(x,y) \in \Gamma$ and $\mu \in (-\mu_0,\mu_0)$,
\begin{align}
\|\mathbf{f}(x,y; \mu)\| & \ge \|\mathbf{f}_0(x,y)\| - |\mu| |c(x,y)| \|\nabla h(x,y)\|
\ge m - |\mu| M > \frac{m}{2} > 0.
\end{align}
Hence,  for every $\mu \in (-\mu_0,\mu_0)$, 
$\Gamma$ is an invariant simple closed curve of~(\ref{eq:DS_0_mu})
containing no equilibria, and is therefore a cycle.

Hyperbolic limit cycles persist and are locally unique under
continuously differentiable perturbations that are sufficiently small 
together with their first derivatives~\cite{Perko}. Therefore,
after decreasing $\mu_0$ if necessary, $\Gamma$ remains hyperbolic 
and preserves its stability for every $\mu \in (-\mu_0,\mu_0)$.
\end{proof}

\section{Algebraic limit cycles in chemical dynamical systems}
\label{sec:4}

In this section, we present some theory and 
results on algebraic limit cycles in general CDSs.
Let us first outline some properties of (semi-)chemical
polynomials, introduced in Definition~\ref{def:semi_chemical}.

\begin{lemma}
\label{lemma:chem_properties}
Let $p, q: \mathbb{R}^2 \to \mathbb{R}$ be polynomials.
Then the following statements hold.
\begin{itemize}
\item[\emph{(i)}] If $p(x,y)$ and $q(x,y)$ are chemical in $x$, 
then, for every $\alpha, \beta \ge 0$, 
$\alpha p(x,y) + \beta q(x,y)$ is chemical in $x$.
\item[\emph{(ii)}] If $p(x,y)$ and $q(x,y)$ are chemical in $x$, 
then $p(x,y) q(x,y)$ is chemical in $x$.
\item[\emph{(iii)}] If $p(x,y)$ is chemical in $x$, 
then the partial derivative 
$\frac{\partial p(x,y)}{\partial y}$ is chemical in $x$.
\item[\emph{(iv)}] For every polynomial $p(x,y)$, 
$x p(x,y)$ is chemical in $x$, 
and $x y p(x,y)$ is chemical in $x$ and $y$.
\end{itemize}
\end{lemma}

\begin{proof}
The results follow directly from Definition~\ref{def:semi_chemical}.
\end{proof}

\noindent
\textbf{Remark}. If $p(x,y)$ is chemical in $x$, 
then $\frac{\partial p(x,y)}{\partial x}$ is 
not necessarily chemical in $x$. In particular,
for every integer $m \ge 0$, $p(x,y) = - x y^m$ is chemical in $x$, 
while $\frac{\partial p(x,y)}{\partial x} = -y^m$ is not.

We now show that any given $n$-degree polynomial can be transformed
into an $(n+1)$-degree (respectively, $(n+2)$-degree)
semi-chemical (respectively, chemical) polynomial 
with all its regular ovals preserved and contained in 
the positive quadrant by first applying a sufficiently large
translation of the coordinates and then multiplying the 
resulting polynomial by $x$ or $y$ (respectively, by $x y$).
We use the notation $\Gamma + (a,b) = \{(x + a,y + b) \in \mathbb{R}^2 | (x,y) \in \Gamma\}$.

\begin{lemma}
\label{lemma:h_construct}
Let $h : \mathbb{R}^2 \to \mathbb{R}$
be a polynomial of degree $n$, and let
$\Gamma \subset \{h = 0\}$ 
be a regular simple closed curve. 
Then, for every sufficiently large 
$a, b > 0$, the polynomial
\begin{align}
H(x,y) & = x^{m} y^{k} h(x - a,y - b)
\label{eq:xy_h}
\end{align}
of degree $(n+m+k)$ satisfies 
$(\Gamma + (a,b)) \subset \{H = 0\} \cap \mathbb{R}_{>}^2$.
Furthermore, $(\ref{eq:xy_h})$ 
is semi-chemical if $m \ge 1$ or $k \ge 1$,
and chemical if $m \ge 1$ and $k \ge 1$.
\end{lemma}

\begin{proof}
By definition, for every $(X,Y) \in \Gamma + (a,b)$,
there exists $(x,y) \in \Gamma$ such that 
$(X,Y) = (x,y) + (a,b)$. Hence, since $h(x,y) = 0$,
\begin{align}
H(X,Y)
& = X^m Y^k 
h(X - a, Y - b) = X^m Y^k  h(x,y) = 0.
\end{align}
Since $\Gamma$ is compact, for every sufficiently large $a,b > 0$, 
its translation $\Gamma + (a,b)$ is in $\mathbb{R}_{>}^2$. 
Further, since $\nabla h(x,y) \ne \mathbf{0}$ on $\Gamma$, 
on the translated curve,  
\begin{align}
\nabla H(X,Y)
& = \nabla (X^m Y^k ) h(x,y)
+ X^m Y^k \nabla h(x,y) = 
X^m Y^k \nabla h(x,y) \ne \mathbf{0}.
\end{align}
Hence, $\Gamma + (a,b)$ is a regular connected component of $\{H = 0\}$. 
The stated chemicality properties follow directly from Lemma~\ref{lemma:chem_properties}(iv).  
\end{proof}

\noindent 
Later in this section, we show that both degree bounds 
can be reduced by one by applying a sufficiently small 
polynomial perturbation and a more general affine change of coordinates;
see Lemmas~\ref{lemma:semichemical_construct} 
and~\ref{lemma:chemical_construct}. 

\textbf{Harnack M-curves}. 
Let $h : \mathbb{R}^2 \to \mathbb{R}$ be a polynomial
of degree $n \ge 1$ regular on $\mathbb{R}^2$, i.e.
$\nabla h(x,y) \ne \mathbf{0}$ for every $(x,y) \in \{h = 0\}$. 
Assume that 
$\{h = 0\}$ contains $M(n)$ 
bounded connected components, where $M(n)$
is given by~(\ref{eq:M_n}).
Then $\{h = 0\}$ is called a \emph{Harnack M-curve} of degree $n$.
For a given degree $n$, Harnack M-curves contain 
the maximum number of bounded 
connected components~\cite{Harnack,Petrovsky}.
We now define chemical analogues of such curves.

\begin{definition}
\label{def:chemical_Harnack}
$($\textbf{Semi-chemical and chemical Harnack M-curves}$)$ 
Let $\{h = 0\}$ be a Harnack M-curve of degree $n$.
If every bounded connected component of $\{h = 0\}$
is contained in the positive quadrant $\mathbb{R}_{>}^2$
and $h$ is semi-chemical $($respectively, chemical$)$, 
then $\{h = 0\}$ is called a \emph{semi-chemical Harnack M-curve} 
$($respectively, a \emph{chemical Harnack M-curve}$)$.
\end{definition} 

\subsection{Dynamics induced by semi-chemical polynomials}
We now use the theory from Section~\ref{sec:3} 
to construct CDSs of the form~(\ref{eq:DS_1}) with semi-chemical 
polynomials $h$ and with prescribed algebraic cycles.

\begin{lemma}
\label{lemma:semi_chemical}
Let $h : \mathbb{R}^2 \to \mathbb{R}$
be an $n$-degree polynomial chemical in $x$, 
and assume that $\{h = 0\}$ contains $N$ 
regular simple closed curves
$\Gamma_1, \Gamma_2, \ldots, 
\Gamma_N  \subset \mathbb{R}_{>}^2$.
Assume that $\alpha_{0,0}, \alpha_{0,1},
\gamma_{0,1}, \gamma_{0,2} \ge 0$
and $\alpha_{1,0}, \beta_{0,1}, \gamma_{1,1} \in \mathbb{R}$ 
are such that, for every $(x,y) \in \mathbb{R}_{>}^2$,
\begin{align}
\gamma_{0,1} + \gamma_{0,2} y - \gamma_{1,1} x & \ne 0 
\; \; \text{and} \; \; 
\left(\alpha_{1,0} \gamma_{0,1} 
+ \alpha_{0,0} \gamma_{1,1} \right)
+ \left( (\alpha_{1,0} - \beta_{0,1}) \gamma_{0,2}
+ \alpha_{0,1} \gamma_{1,1} \right) y \ne 0.
\label{eq:conditions_semichemical}
\end{align}
Then system
\begin{align}
\frac{\mathrm{d} x}{\mathrm{d} t} 
& = (\alpha_{0,0} + \alpha_{1,0} x + \alpha_{0,1} y) h(x,y) - 
\left(-\gamma_{0,1} y + \gamma_{1,1} x y - \gamma_{0,2} y^2 \right) 
\frac{\partial h(x,y)}{\partial y}, 
\nonumber \\
\frac{\mathrm{d} y}{\mathrm{d} t} 
& = (\beta_{0,1} y) h(x,y) 
+ \left(-\gamma_{0,1} y + \gamma_{1,1} x y - \gamma_{0,2} y^2 \right) 
\frac{\partial h(x,y)}{\partial x},
\label{eq:DS_h_semichemical}
\end{align}  
is a \emph{CDS} of degree at most $n+1$ and 
admits $N$ hyperbolic algebraic limit cycles 
$\Gamma_1, \Gamma_2, \ldots, \Gamma_N$.
For each $i = 1,2,\ldots, N$, the cycle $\Gamma_i$
is stable (respectively, unstable) if 
$K_{\Gamma_i} < 0$ (respectively, $K_{\Gamma_i} > 0$),
where $K_{\Gamma_i}$ is given by~$(\ref{eq:K_F})$.
Furthermore, $(\ref{eq:DS_h_semichemical})$
has no other cycles in $\mathbb{R}_{>}^2$
if $\{h = 0\}$ contains no 
bounded connected components in $\mathbb{R}_{>}^2$
other than $\Gamma_1, \Gamma_2, \ldots, \Gamma_N$.
\end{lemma}

\begin{proof}
By Lemma~\ref{lemma:chem_properties},
(\ref{eq:DS_h_semichemical}) is a CDS
of degree at most $n+1$. 
It is of the form~(\ref{eq:DS_1}),
with $f_1(x,y) = (\alpha_{0,0} + \alpha_{1,0} x + \alpha_{0,1} y)$, 
$g_1(x,y) = \beta_{0,1} y$ and 
$c(x,y) = (-\gamma_{0,1} - \gamma_{0,2} y + \gamma_{1,1} x) y$.
Let $F$ be as defined in~(\ref{eq:F}).
By~(\ref{eq:conditions_semichemical}),
$c(x,y) \ne 0$ and $F(x,y) \ne 0$
for every $(x,y) \in \mathbb{R}_{>}^2$.
Hence, the conclusions follow from
Theorem~\ref{theorem:AL_strong} with 
$\mathbb{S} = \mathbb{R}_{>}^2$.
\end{proof}

\noindent
\textbf{Remark}. 
Conditions~(\ref{eq:conditions_semichemical}) 
are satisfied if
$\alpha_{0,0} = \gamma_{1,1} = 1$ 
and $\alpha_{1,0} = \alpha_{0,1} 
= \beta_{0,1} = \gamma_{0,1} = \gamma_{0,2} = 0$,
with~(\ref{eq:DS_h_semichemical}) given by
\begin{align}
\frac{\mathrm{d} x}{\mathrm{d} t} 
& = h(x,y) - x y \frac{\partial h(x,y)}{\partial y}, \nonumber \\
\frac{\mathrm{d} y}{\mathrm{d} t} 
& = x y \frac{\partial h(x,y)}{\partial x}.
\label{eq:DS_h_semichemical_1}
\end{align}

Let us characterize when the $(n+1)$-degree terms cancel out 
in CDS~(\ref{eq:DS_h_semichemical}).

\begin{lemma}
$($\textbf{Cancellation lemma}$)$ 
\label{lemma:cancellation}
Let $h_n$ denote the $n$-degree homogeneous part of polynomial $h$.
Under the assumptions of \emph{Lemma~\ref{lemma:semi_chemical}}, 
\emph{CDS}~$(\ref{eq:DS_h_semichemical})$ has degree $n$
if and only if $h_n$ is a single monomial,
\begin{align}
h_n(x,y) & = \kappa_{i,j} x^i y^j, \; \; \text{where} 
\; \; \kappa_{i,j} \ne 0 \; \; \text{and} \; \; i + j = n, 
\label{eq:semi_chemical_cancel_1}
\end{align}
and the coefficients in~$(\ref{eq:DS_h_semichemical})$ are chosen to satisfy
\begin{align}
\alpha_{0,1} & = \gamma_{0,2} = 0, \; \; \; 
\alpha_{1,0} = j \gamma_{1,1}, \; \; \; 
\beta_{0,1} = - i \gamma_{1,1}.
\label{eq:semi_chemical_cancel_2}
\end{align}
Otherwise, \emph{CDS}~$(\ref{eq:DS_h_semichemical})$ has degree $n+1$. 
\end{lemma}

\begin{proof}
Let $f_{n+1}$ and $g_{n+1}$ denote the $(n+1)$-degree
homogeneous parts of the two right-hand 
sides in~(\ref{eq:DS_h_semichemical}). Then
\begin{align}
f_{n+1}(x,y) & = (\alpha_{1,0} x + \alpha_{0,1} y) h_n(x,y) - 
\left(\gamma_{1,1} x y - \gamma_{0,2} y^2 \right) 
\frac{\partial h_n(x,y)}{\partial y}, \nonumber \\
g_{n+1}(x,y) & = (\beta_{0,1} y) h_n(x,y) 
+ \left(\gamma_{1,1} x y - \gamma_{0,2} y^2 \right) 
\frac{\partial h_n(x,y)}{\partial x}. \label{eq:fg_n1}
\end{align}
By homogeneity, $h_n(\lambda x, \lambda y) = \lambda^n h_n(x,y)$
for every $\lambda > 0$, 
it follows that $x \frac{\partial h_n}{\partial x}
+ y \frac{\partial h_n}{\partial y} = n h_n$. 
Consequently, 
\begin{align}
x g_{n+1}(x,y) - y f_{n+1}(x,y) & = y \big[ 
(\beta_{0,1} - \alpha_{1,0} + n \gamma_{1,1}) x
- (\alpha_{0,1} + n \gamma_{0,2}) y \big] h_n(x,y).
\label{eq:cancellation_identity_1}
\end{align}  

Assume that the $(n+1)$-degree terms cancel out, $f_{n+1} = g_{n+1} = 0$.
Then, by~(\ref{eq:cancellation_identity_1}), since $h_n \not \equiv 0$, 
$\alpha_{1,0} = \beta_{0,1} + n \gamma_{1,1}$ and 
$\alpha_{0,1} + n \gamma_{0,2} = 0$. By
assumption, $\alpha_{0,1}, \gamma_{0,2} \ge 0$, 
so the latter constraint gives $\alpha_{0,1} = \gamma_{0,2} = 0$.
Furthermore, $\gamma_{1,1} \ne 0$ since, otherwise, 
$f_{n+1} = g_{n+1} = 0$ would imply that 
$\alpha_{1,0} = \beta_{0,1} = 0$, contradicting 
conditions~(\ref{eq:conditions_semichemical}).
Consequently, the cancellation conditions reduce to
\begin{align}
0 & = \alpha_{1,0} x h_n(x,y) - 
\gamma_{1,1} x y \frac{\partial h_n(x,y)}{\partial y}, \; \; \; \;
0 = \beta_{0,1} y h_n(x,y) + \gamma_{1,1} x y
\frac{\partial h_n(x,y)}{\partial x}. \nonumber 
\end{align}
Substituting $h_n(x,y) = \sum_{i + j = n} \kappa_{i,j} x^i y^j$ 
and comparing the coefficients, it follows that, 
for every $i + j = n$ with $\kappa_{i,j} \ne 0$, 
$\alpha_{1,0} = j \gamma_{1,1}$ and $\beta_{0,1} = - i \gamma_{1,1}$. 
Since $\gamma_{1,1} \ne 0$, these identities can hold for
at most one pair $(i,j)$, implying~(\ref{eq:semi_chemical_cancel_1})
and~(\ref{eq:semi_chemical_cancel_2}).
Conversely, if conditions~(\ref{eq:semi_chemical_cancel_1})
and~(\ref{eq:semi_chemical_cancel_2}) hold, then the cancellation
of $(n+1)$-degree terms follows directly from~(\ref{eq:fg_n1}).
Furthermore, the coefficient of $x^i y^j$, where $i + j = n$, 
in the first equation of~(\ref{eq:DS_h_semichemical}) is given by 
$(\alpha_{0,0} + j \gamma_{0,1}) \kappa_{i,j} + 
(\alpha_{1,0} - j \gamma_{1,1}) \kappa_{i-1,j} = 
(\alpha_{0,0} + j \gamma_{0,1}) \kappa_{i,j}$.
By~(\ref{eq:conditions_semichemical}),
$\alpha_{0,0} + j \gamma_{0,1} \ne 0$. 
Therefore, CDS~(\ref{eq:DS_h_semichemical}) has degree $n$. 
\end{proof}

\noindent 
\textbf{Remark}. Conditions~(\ref{eq:conditions_semichemical}) 
and~(\ref{eq:semi_chemical_cancel_2})
are satisfied if $\alpha_{1,0} = j$, $\beta_{0,1} = -i$
with $i + j = n$, $\alpha_{0,0} = \gamma_{1,1} = 1$ 
and $\alpha_{0,1} = \gamma_{0,1} = \gamma_{0,2} = 0$,
with~(\ref{eq:DS_h_semichemical}) given by
\begin{align}
\frac{\mathrm{d} x}{\mathrm{d} t} 
& = (1 + j x) h(x,y) - x y  
\frac{\partial h(x,y)}{\partial y}, 
\nonumber \\
\frac{\mathrm{d} y}{\mathrm{d} t} 
& = - i y h(x,y) 
+ x y \frac{\partial h(x,y)}{\partial x},
\label{eq:DS_h_semichemical_2}
\end{align} 
where $i + j = n$. If $h$ satisfies~(\ref{eq:semi_chemical_cancel_1}),
then CDS~(\ref{eq:DS_h_semichemical_2}) has degree $n$
by Lemma~\ref{lemma:cancellation}. 

\noindent 
\textbf{Remark}. For every $n \ge 3$, an invariant algebraic curve
$\{h = 0\}$ satisfying~(\ref{eq:semi_chemical_cancel_1})
is non-generic in the sense of~\cite{Llibre}, because 
its highest-degree homogeneous part is a single monomial
with repeated factors. 

\subsubsection{Construction of semi-chemical polynomials and Harnack M-curves}
We now prove \emph{constructively} that
any given $n$-degree polynomial can be transformed
into an $n$-degree semi-chemical polynomial 
while preserving its regular ovals 
and placing them in the positive quadrant.
In what follows, we measure the distance between two non-empty compact sets 
$A, B \subset \mathbb{R}^2$ using the Hausdorff distance, defined by
$d_H(A,B) = \max \{\max_{\mathbf{x} \in A} \textrm{dist}(\mathbf{x},B), 
\max_{\mathbf{y} \in B} \textrm{dist}(\mathbf{y},A)\}$.

\begin{lemma}
\label{lemma:semichemical_construct}
Let $H : \mathbb{R}^2 \to \mathbb{R}$ 
be a polynomial of degree $n$ such that
$\{H = 0\}$ contains $N \ge 1$ regular simple closed curves.
Then there exists a semi-chemical 
polynomial of degree $n$ 
whose zero set contains $N$ regular 
simple closed curves in $\mathbb{R}_{>}^2$.
\end{lemma}

\begin{proof}
Let $H(u,v)$ be an $n$-degree polynomial such that
$\{H = 0\}$ contains $N \ge 1$ regular simple closed curves
$\Gamma_i$, $i = 1,2,\ldots,N$. 
Let $\gamma_{n,0}$ denote the coefficient of $u^n$ in $H$. 
Without loss of generality, possibly after an arbitrarily small polynomial
perturbation, we may assume that 
\begin{align}
(-1)^n \gamma_{n,0} > 0. \label{eq:sign_gamma}
\end{align}
Indeed, if $(-1)^n \gamma_{n,0} < 0$, then replace
$H$ with $-H$. Otherwise, if $\gamma_{n,0} = 0$, consider the perturbed 
polynomial $H_{\varepsilon}(u,v) = H(u,v) + \varepsilon u^n$
with $\varepsilon > 0$. 
Since the curves $\Gamma_i$, $i = 1,2,\ldots,N$,
are regular and pairwise disjoint, the implicit function theorem and compactness
imply that, for every sufficiently small $\varepsilon > 0$, 
$\{H_{\varepsilon} = 0\}$ contains $N$ pairwise disjoint
regular simple closed curves $\Gamma_{i,\varepsilon}$ satisfying 
$d_H(\Gamma_{i,\varepsilon},\Gamma_i) \to 0$ as $\varepsilon \to 0$
for every $i = 1,2,\ldots,N$.
If $(-1)^n \varepsilon > 0$, then replace $H$ with $H_{\varepsilon}$;
otherwise, replace it with $-H_{\varepsilon}$.
(Alternatively, the condition $\gamma_{n,0} \ne 0$ can be achieved
via a suitable rotation of the coordinates.) 

Let $\mathbb{K} \subset \mathbb{R}^2$ be a compact set containing all 
bounded connected components of $\{H = 0\}$. Choose $b > 0$
such that $v + b > 0$ for every $(u,v) \in \mathbb{K}$. 
Fix this $b$ and, for every $a > 0$, define
\begin{align}
h(x,y) & = H(x - a (1 + y), y - b) = \sum_{i + j = 0}^n \gamma_{i,j} 
\left(x - a (1 + y) \right)^i (y - b)^j.
\end{align}
At $x = 0$, applying the binomial theorem, one obtains
\begin{align}
h(0,y) 
& =  \sum_{i + j = 0}^n \gamma_{i,j} (-a)^i (1 + y)^i (y - b)^j
 = \sum_{i + j = 0}^n \sum_{r = 0}^i \sum_{l = 0}^j 
\gamma_{i,j} (-a)^i (-b)^{j-l} \binom{i}{r} \binom{j}{l} y^{r+l} \nonumber \\
& = (-1)^n \gamma_{n,0} a^n \sum_{k = 0}^n \binom{n}{k} y^{k}
+ \sum_{\substack{i + j = 0, \\ (i,j) \ne (n,0)}}^n \sum_{r = 0}^i \sum_{l = 0}^j 
\gamma_{i,j} (-a)^i (-b)^{j-l} \binom{i}{r} \binom{j}{l} y^{r+l}.
\end{align}
Hence, for every $k = 0,1,\ldots,n$, the coefficient of $y^k$
is a polynomial in $a$ of the form 
\begin{align}
(-1)^n \gamma_{n,0} \binom{n}{k} a^n + R_k(a),
\end{align}
where $R_k(a)$ is zero or has degree at most $n-1$. Since 
$(-1)^n \gamma_{n,0} > 0$, for every sufficiently large
$a > 0$, all coefficients of $h(0,y)$ are positive,
so that $h(x,y)$ is chemical in $x$ by Definition~\ref{def:semi_chemical}. 

The affine change of coordinates $x = u + a (1 + v + b)$, 
$y = v + b$ is invertible. On $\mathbb{K}$, $y > 0$
and $(1 + v + b) > 1$. Hence, by further increasing $a$ if necessary, 
$x > 0$ on $\mathbb{K}$, so that all bounded connected
components of $\{H = 0\}$ are mapped into $\mathbb{R}_{>}^2$.
Finally, since an invertible affine 
change of coordinates preserves the degree, regularity
and the number of bounded connected components,
$h(x,y)$ is an $n$-degree semi-chemical polynomial 
whose zero set contains at least $N$ 
regular simple closed curves in $\mathbb{R}_{>}^2$.
\end{proof}

Lemma~\ref{lemma:semichemical_construct} implies the following result.

\begin{theorem}
\label{theorem:semi_chemical_Harnack}
For every $n \ge 2$, there exists a semi-chemical Harnack M-curve of degree $n$. 
\end{theorem}

\begin{proof}
Let $H(u,v)$ be a polynomial of degree $n$ such that $\{H = 0\}$ is a Harnack M-curve
with bounded connected components $\Gamma_i$, $i = 1,2,\ldots,M(n)$.
The conclusion follows by applying the construction from 
Lemma~\ref{lemma:semichemical_construct},
and noting that a Harnack M-curve
of degree $n$ remains a Harnack M-curve under sufficiently small 
polynomial perturbations of degree at most $n$~\cite{Viro}.  
\end{proof}

Lemma~\ref{lemma:semi_chemical} and Theorem~\ref{theorem:semi_chemical_Harnack}
together give the following result.

\begin{corollary}
\label{corollary:semi_chemical_n}
For every integer $n \ge 2$, $H_a^c(n) \ge M(n-1)$,
where $M(n)$ is given by~$(\ref{eq:M_n})$.
\end{corollary}

\begin{proof}
For $n = 2$, the inequality follows from $M(1) = 0$. Assume that $n \ge 3$.
Let $h$ be an $(n-1)$-degree polynomial such that $\{h = 0\}$
is a semi-chemical Harnack M-curve, which exists 
by Theorem~\ref{theorem:semi_chemical_Harnack}. 
Interchanging coordinates if necessary,
assume that $h$ is chemical in $x$.
By Lemma~\ref{lemma:cancellation}, 
CDS~(\ref{eq:DS_h_semichemical_1}) has degree $n$,
since $\alpha_{1,0} = \beta_{0,1} = 0$.
Furthermore, by Lemma~\ref{lemma:semi_chemical}, 
this CDS admits $M(n - 1)$ algebraic limit cycles,
which are the ovals of $\{h = 0\}$ in the positive quadrant.
Hence $H_a^c(n) \ge M(n-1)$.
\end{proof}

\subsubsection{Semi-chemical polynomials with a single leading monomial}
We now strengthen Corollary~\ref{corollary:semi_chemical_n} 
by proving the existence of suitable semi-chemical polynomials  
with a single leading monomial, which, using Lemma~\ref{lemma:cancellation}, 
allow the cancellation of degree $n+1$ terms in~(\ref{eq:DS_h_semichemical}).

\begin{theorem}
\label{theorem:semi_chemical_Harnack_stronger}
For every $n \ge 3$, there exists an $n$-degree 
semi-chemical polynomial $h : \mathbb{R}^2 \to \mathbb{R}$ 
with a single leading monomial, and 
such that the zero set $\{h = 0\}$ is regular and contains exactly
\begin{align}
\begin{cases}
M(n) - 1, & \text{if } n \text{ is even},  \\
M(n), & \text{if } n \text{ is odd}, 
\end{cases} 
\label{eq:semi_chemical_improved}
\end{align} 
ovals in $\mathbb{R}_{>}^2$, where $M(n)$ is given by~$(\ref{eq:M_n})$.
In particular, for every odd $n \ge 3$, there exists a semi-chemical 
Harnack M-curve of degree $n$ with a single leading monomial. 
\end{theorem}

\begin{proof}
Let $n \ge 3$, and let $\Delta_n$ be the triangle with vertices
$(0,0)$, $(n-1,0)$ and $(0,n)$. Its area is $A_n = n (n-1)/2$.
Furthermore, its horizontal and vertical sides contain $n$ and $(n+1)$
integer lattice points, respectively, while this number 
for the diagonal side is $2$, giving in total 
$B_n = n + (n+1) + 2 - 3 = 2 n$ lattice points
on the boundary of $\Delta_n$. By Pick's theorem, the total number
of interior integer lattice points in $\Delta_n$ is then given by 
$I_n = A_n - B_n/2 + 1 = (n-1) (n-2)/2$, 
which is equal to~(\ref{eq:semi_chemical_improved}). 
By~\cite{Mikhalkin}[Corollaries A4, 9 and 10], 
there exists a real $n$-degree polynomial 
\begin{align}
H(u,v) & = \sum_{i+j = 0}^n \gamma_{i,j} u^i v^j = 
\gamma_{0,n} v^n + \gamma_{n-1,0} u^{n-1} + R(u,v),
\label{eq:H_improved}
\end{align}
with $\gamma_{0,0}, \gamma_{n-1,0}, \gamma_{0,n} \ne 0$ and 
$\gamma_{i,j} = 0$ for every $(i,j) \not \in \Delta_n \cap \mathbb{Z}_{\ge}^2$,
and such that $\{H = 0\}$ is regular and contains at least $I_n$ ovals.
Note that $R(u,v)$ is of degree at most $n-1$ and the exponent 
of $u$ in each of its monomials is at most $n-2$.

We now show that $\{H=0\}$ contains exactly $I_n$ ovals. 
To this end, consider the projective curve defined by the homogenization
$\mathcal{H}(U,V,W)=W^n H(U/W,V/W)$. 
Since the only non-zero monomial of degree $n$ in
$H$ is $\gamma_{0,n}v^n$, we have
$\mathcal{H}(U,V,0)=\gamma_{0,n}V^n.$
Hence, the projective curve has a unique point at
infinity, namely $P=[1:0:0]$. Since 
$\frac{\partial\mathcal{H}}{\partial W}(1,0,0)
=\gamma_{n-1,0}\neq 0$, $P$ is a regular point.
Consequently, the real projective curve is regular
and, by Harnack's theorem, has at most 
$(n-1)(n-2)/2 + 1 = I_n + 1$
connected components. If necessary, this bound
can be obtained by considering a sufficiently small
generic perturbation of $\mathcal{H}$ that preserves
the topology of its real zero set and yields a
non-singular complex projective curve.
Because one of the connected components 
of the original projective curve contains $P$, 
and therefore does not correspond to an oval
of the affine curve $\{H=0\}$, the latter contains
at most $I_n$ ovals. Since $\{H=0\}$ contains 
at least $I_n$ ovals by construction, 
it contains exactly $I_n$ ovals.

Without loss of generality, we assume that 
the coefficients in~(\ref{eq:H_improved}) satisfy 
\begin{align}
\gamma_{0,n} & > 0, \; \; \; \; \; \; (-1)^{n-1} \gamma_{n-1,0} > 0.
\end{align}
Indeed, assume first that $n$ is odd. 
If $(-1)^{n-1} \gamma_{n-1,0} < 0$, then consider $- H(u,v)$ instead.
If $\gamma_{0,n} < 0$, then consider instead $H(u,-v)$,
which ensures that $\gamma_{0,n} > 0$ without changing $(-1)^{n-1} \gamma_{n-1,0}$.
Assume now that $n$ is even. 
If $\gamma_{0,n} < 0$, then consider $- H(u,v)$ instead.
If $(-1)^{n-1} \gamma_{n-1,0} < 0$, then consider instead $H(-u,v)$, 
which ensures that $(-1)^{n-1} \gamma_{n-1,0} > 0$ without changing $\gamma_{0,n}$.

We now proceed as in the proof of Theorem~\ref{theorem:semi_chemical_Harnack}. 
Let $\mathbb{K} \subset \mathbb{R}^2$ be a compact set containing all 
the ovals of $\{H = 0\}$. Choose $b > 0$
such that $v + b > 0$ for every $(u,v) \in \mathbb{K}$. 
Fix this $b$ and, for every $a > 0$, define
$h(x,y) = H(x - a (1 + y), y - b)$.
For every sufficiently large $a$, all the ovals
of $\{h = 0\}$ are contained in $\mathbb{R}_{>}^2$. 
Setting $x = 0$ and using the binomial theorem, one obtains 
\begin{align}
h(0,y) & = \gamma_{0,n} (y - b)^n + (-1)^{n-1} \gamma_{n-1,0} a^{n-1} (1 + y)^{n-1} 
+ R(- a (1 + y),y - b) \nonumber \\
& = \gamma_{0,n} \sum_{k = 0}^n \binom{n}{k} (-b)^{n-k} y^k
+ (-1)^{n-1} \gamma_{n-1,0} a^{n-1} \sum_{k = 0}^{n-1}  \binom{n-1}{k} y^k
+ \sum_{k = 0}^{n-1} R_k(a) y^k \nonumber \\
& = \gamma_{0,n} y^n + \sum_{k = 0}^{n-1} 
\Big[\binom{n}{k} \gamma_{0,n} (-b)^{n-k}
+ \binom{n-1}{k} (-1)^{n-1} \gamma_{n-1,0} a^{n-1} + R_k(a) \Big] y^k.
\end{align}
Since $R$ is of degree at most $n-2$ in its first argument, 
each $R_k(a)$ is a polynomial in $a$ of degree at most $n-2$.
Hence, since $\gamma_{0,n} > 0$ and $(-1)^{n-1} \gamma_{n-1,0} > 0$, 
it follows that, for every sufficiently large $a > 0$, 
all coefficients of $h(0,y)$ are positive. Furthermore, the only non-zero
monomial of degree $n$ in $h$ is $\gamma_{0,n} y^n$. 
Therefore, $h$ is semi-chemical, with a single leading monomial,
and with the number of ovals in $\mathbb{R}_{>}^2$ 
given by~(\ref{eq:semi_chemical_improved}). 
\end{proof}

\begin{theorem}
\label{theorem:semi_chemical_improved}
For every integer $n \ge 2$, $H_a^c(n) \ge M(n) - (1 + (-1)^n)/2$,
where $M(n)$ is given by~$(\ref{eq:M_n})$.
\end{theorem}

\begin{proof}
For $n = 2$, the inequality follows from $M(2) - 1 = 0$. Assume that $n \ge 3$.
Let $h$ be an $n$-degree semi-chemical polynomial
with a single leading monomial such that $\{h = 0\}$ is regular
and contains the number of ovals in $\mathbb{R}_{>}^2$ 
given by~(\ref{eq:semi_chemical_improved}),
which exists by Theorem~\ref{theorem:semi_chemical_Harnack_stronger}.
Interchanging coordinates if necessary,
assume that $h$ is chemical in $x$.
By Lemma~\ref{lemma:cancellation}, 
CDS~(\ref{eq:DS_h_semichemical_2}) has degree $n$,
and, by Lemma~\ref{lemma:semi_chemical}, 
it admits $M(n) - 1$ hyperbolic limit cycles in $\mathbb{R}_{>}^2$ 
if $n$ is even, and $M(n)$ if $n$ is odd, proving the theorem.
\end{proof}

\subsection{Dynamics induced by chemical polynomials}
Let us now present results involving chemical polynomials.

\begin{lemma}
\label{lemma:chemical}
Let $h : \mathbb{R}^2 \to \mathbb{R}$
be an $n$-degree chemical polynomial and assume that
$\{h = 0\}$ contains $N$ regular simple closed curves
$\Gamma_1, \Gamma_2, \ldots, 
\Gamma_N  \subset \mathbb{R}_{>}^2$.
Assume that  $\alpha_{0,0}, \beta_{0,0},
\gamma_{1,0}, \gamma_{0,1} \ge 0$,
\begin{align}
\gamma_{1,0} \gamma_{0,1} & = 0 \; \; \text{and} \; \; 
(\alpha_{0,0} \gamma_{1,0} - \beta_{0,0} \gamma_{0,1}) \ne 0.
\label{eq:conditions_chemical}
\end{align}
Then system
\begin{align}
\frac{\mathrm{d} x}{\mathrm{d} t} 
& = \alpha_{0,0} h(x,y) - 
\left(\gamma_{1,0} x - \gamma_{0,1} y \right)
 \frac{\partial h(x,y)}{\partial y}, 
\nonumber \\
\frac{\mathrm{d} y}{\mathrm{d} t} 
& = \beta_{0,0} h(x,y) 
+ \left(\gamma_{1,0} x - \gamma_{0,1} y \right) 
\frac{\partial h(x,y)}{\partial x},
\label{eq:DS_h_chemical}
\end{align} 
is a \emph{CDS} of degree $n$ and admits 
$N$ hyperbolic algebraic limit cycles 
$\Gamma_1, \Gamma_2, \ldots, \Gamma_N$.
For each $i = 1,2,\ldots, N$, the cycle $\Gamma_i$
is stable (respectively, unstable) if 
$h(\mathbf{q}_i) > 0$
(respectively, $h(\mathbf{q}_i) < 0$),
where $\mathbf{q}_i$ is any point
sufficiently close to $\Gamma_i$ and in 
the exterior of $\Gamma_i$.
Furthermore, $(\ref{eq:DS_h_chemical})$
has no other cycles in $\mathbb{R}_{>}^2$
if $\{h = 0\}$ contains no 
bounded connected components in $\mathbb{R}_{>}^2$
other than $\Gamma_1, \Gamma_2, \ldots, \Gamma_N$.
\end{lemma}

\begin{proof}
By Lemma~\ref{lemma:chem_properties}, (\ref{eq:DS_h_chemical}) is a CDS. 
Let $h_n$ denote the $n$-degree homogeneous part of $h$. 
By~(\ref{eq:conditions_chemical}), either $\gamma_{1,0} = 0$ or $\gamma_{0,1} = 0$.
Assume first that $\gamma_{1,0} = 0$, so that $\beta_{0,0}, \gamma_{0,1} > 0$.
The $n$-degree homogeneous part of the second equation from~(\ref{eq:DS_h_chemical})
is then $\beta_{0,0} h_n - \gamma_{0,1} y \frac{\partial h_n}{\partial x}$.
Let $\kappa_{i,n-i} x^i y^{n-i}$ be the monomial of $h_n$ 
with $\kappa_{i,n-i} \ne 0$ and the largest index $i \ge 0$. Since every monomial of 
$y \frac{\partial h_n}{\partial x}$ has $x$-degree at most $i - 1$, 
it follows that the coefficient of $x^i y^{n-i}$ is 
given by $\beta_{0,0} \kappa_{i,n-i} \ne 0$. 
Hence, the second equation is of degree $n$. 
Analogous argument shows that, if $\gamma_{0,1} = 0$,
then the first equation is of degree $n$.
Therefore, (\ref{eq:DS_h_chemical}) is of degree $n$
under the conditions~(\ref{eq:conditions_chemical}).

CDS(\ref{eq:DS_h_chemical}) is of the form~(\ref{eq:DS_1}),
with $f_1(x,y) = \alpha_{0,0}$, $g_1(x,y) =  \beta_{0,0}$ 
and $c(x,y) = \left(\gamma_{1,0} x - \gamma_{0,1} y \right)$.
Let $F$ be as defined in~(\ref{eq:F}).
By~(\ref{eq:conditions_chemical}), $c(x,y) \ne 0$ and 
$F(x,y) = - (\alpha_{0,0} \gamma_{1,0} - \beta_{0,0} \gamma_{0,1}) \ne 0$ 
for every $(x,y) \in \mathbb{R}_{>}^2$.
Furthermore, conditions~(\ref{eq:conditions_chemical}) imply
that $\textrm{sign}(F) = - \textrm{sign}(c)$ on $\mathbb{R}_{>}^2$. 
Hence, by~(\ref{eq:K_F}), $\textrm{sign}(K_{\Gamma_i}) 
= s(\mathbf{x}_{\Gamma_i}) \textrm{sign}(F)
= - \textrm{sign}(h(\mathbf{q}_i))$,
where $\mathbf{q}_i$ is any point immediately outside $\Gamma_i$. 
Hence, the conclusions follow from
Theorem~\ref{theorem:AL_strong} with 
$\mathbb{S} = \mathbb{R}_{>}^2$.
\end{proof}

\noindent
\textbf{Remark}. 
Conditions~(\ref{eq:conditions_chemical}) 
are satisfied if
$\alpha_{0,0} = \gamma_{1,0} = 0$ 
and $\beta_{0,0} = \gamma_{0,1} = 1$,
with~(\ref{eq:DS_h_chemical}) given by
\begin{align}
\frac{\mathrm{d} x}{\mathrm{d} t} 
& =  y \frac{\partial h(x,y)}{\partial y}, 
\nonumber \\
\frac{\mathrm{d} y}{\mathrm{d} t} 
& = h(x,y) - y  \frac{\partial h(x,y)}{\partial x}.
\label{eq:DS_h_chemical_1}
\end{align}

\subsubsection{Construction of chemical polynomials}
Applying Lemma~\ref{lemma:semichemical_construct}, 
one obtains the following result.

\begin{lemma}
\label{lemma:chemical_construct}
Let $H : \mathbb{R}^2 \to \mathbb{R}$ 
be a polynomial of degree $n$ such that
$\{H = 0\}$ contains $N \ge 1$ regular simple closed curves.
Then there exists a chemical 
polynomial of degree $n+1$ 
whose zero set contains $N$ regular 
simple closed curves in $\mathbb{R}_{>}^2$.
\end{lemma}

\begin{proof}
By Lemma~\ref{lemma:semichemical_construct}, 
there exists an $n$-degree polynomial $h(x,y)$
chemical in $x$ such that $\{h = 0\}$ contains 
$N$ regular simple closed curves in $\mathbb{R}_{>}^2$.
Define the $(n+1)$-degree polynomial $\tilde{h}(x,y) = y h(x,y)$.
By Lemma~\ref{lemma:chem_properties}(iv),
$\tilde{h}$ is chemical. Furthermore, since $y > 0$
in $\mathbb{R}_{>}^2$ and $\nabla \tilde{h} = y \nabla h$ on $\{h = 0\}$, 
multiplication by $y$ preserves 
the $N$ regular simple closed curves of $\{h = 0\}$ 
in $\mathbb{R}_{>}^2$, completing the proof.
\end{proof}

\noindent
\textbf{Remark}. Corollary~\ref{corollary:semi_chemical_n}
can also be derived from 
Lemmas~\ref{lemma:chemical} and~\ref{lemma:chemical_construct}.

Lemma~\ref{lemma:chemical_construct} cannot be used directly
to construct chemical Harnack M-curves, since it
increases the polynomial degree by one.
We now investigate the existence 
of chemical Harnack M-curves of degrees $n = 2,3,4$.

\begin{lemma}
\label{lemma:cHarnack_2}
There exists no chemical Harnack M-curve of degree $n = 2$. 
\end{lemma}

\begin{proof}
Assume for contradiction that 
$\{h = 0\}$ is a chemical Harnack M-curve, where
$h(x,y) = \kappa_{0,0} + \kappa_{1,0} x + \kappa_{0,1} y
+ \kappa_{2,0} x^2 + \kappa_{1,1} x y + \kappa_{0,2} y^2$. 
By chemicality, $\kappa_{0,0}, \kappa_{1,0}, \kappa_{0,1}, 
\kappa_{2,0}, \kappa_{0,2} \ge 0$. 
A quadratic Harnack M-curve consists of a single ellipse. 
Consequently, its quadratic part
$Q(x,y) = \kappa_{2,0} x^2 + \kappa_{1,1} x y + \kappa_{0,2} y^2$
is definite. Since $\kappa_{2,0}, \kappa_{0,2} \ge 0$, $Q$
is positive definite, i.e. $Q(x,y) > 0$ for every $(x,y) \ne (0,0)$. 
Together with chemicality, it follows that $h(x,y) > 0$ 
for every $(x,y) \in \mathbb{R}_{>}^2$. Therefore $\{h = 0\}$ 
does not intersect the positive quadrant, a contradiction.
\end{proof}

\begin{lemma}
\label{lemma:cHarnack_3}
There exists a chemical Harnack M-curve of degree $n = 3$. 
\end{lemma}

\begin{proof}
Define the cubic polynomial
\begin{align}
H(x,y) & = 1 - 3 x y + x^3 + y^3, \nonumber 
\end{align}
and, for every $\varepsilon > 0$, define 
\begin{align}
h_{\varepsilon}(x,y) = H(x,y) - \varepsilon.
\nonumber 
\end{align}
For every $0 < \varepsilon < 1$, 
$h_{\varepsilon}$ is chemical by 
Definition~\ref{def:chemical_Harnack} because
\begin{align}
h_{\varepsilon}(0,y) =  (1 - \varepsilon) + y^3,
\; \; \; \; 
h_{\varepsilon}(x,0) & = (1 - \varepsilon) + x^3.
\nonumber 
\end{align}

The points satisfying $\nabla H = \mathbf{0}$
are given by $(0,0)$ and $(1,1)$,
and $H(1,1) = 0$. Since the
Hessian matrix $\nabla^2 H(1,1)$ has positive eigenvalues, 
$(1,1)$ is a non-degenerate minimum of $H$. 

Consider the auxiliary Hamiltonian DS
\begin{align}
\frac{\mathrm{d} x}{\mathrm{d} t} & = 
\frac{\partial H(x,y)}{\partial y}, \; \; \; \; \; \; 
\frac{\mathrm{d} y}{\mathrm{d} t} = 
- \frac{\partial H(x,y)}{\partial x}.
\label{eq:auxiliary_H}
\end{align}
Since $(1,1)$ is a non-degenerate minimum of $H$, 
it is a non-degenerate center of this DS~\cite{Perko}. 
Therefore, for every sufficiently small $\varepsilon > 0$, 
the level set $\{H = \varepsilon\}$ contains a
cycle $\Gamma \subset \mathbb{R}_{>}^2$ 
surrounding the center. 
This cycle is a regular oval of $\{h_{\varepsilon} = 0\}$
because $\{H = \varepsilon\} = \{h_{\varepsilon} = 0\}$.
Furthermore, the two points satisfying
$\nabla H = \mathbf{0}$ satisfy $H = 0$ or $H = 1$.
Hence, for every $0 < \varepsilon < 1$, 
$\nabla H \ne \mathbf{0}$ on $\{H = \varepsilon\}$,
i.e. $\{h_{\varepsilon} = 0\}$ is regular.
Since a regular cubic has at most $M(3) = 1$ oval,
$\Gamma$ is the only bounded connected component
of $\{h_{\varepsilon} = 0\}$, and the Harnack bound is attained.
Therefore $\{h_{\varepsilon} = 0\}$ is a chemical Harnack M-curve.
\end{proof}

\begin{lemma}
\label{lemma:cHarnack_4}
There exists a chemical Harnack M-curve of degree $n = 4$. 
\end{lemma}

\begin{proof}
Define the quadratic polynomials
\begin{align}
A(x,y) & = (x + y - 2) (x + y - 3), \nonumber \\
B(x,y) & = 2 (x - y)^2 - (x + y)^2 + x + y + 1,
\nonumber 
\end{align}
and, for every $\varepsilon > 0$, define the associated quartic polynomials
\begin{align}
H(x,y) & = A(x,y)^2 + 31 B(x,y)^2, \; \; \; \; \; \; 
h_{\varepsilon}(x,y) = H(x,y) - \varepsilon.
\label{eq:h_quartic_single}
\end{align}
For every $0 < \varepsilon < 67$,  
$h_{\varepsilon}$ is chemical by 
Definition~\ref{def:chemical_Harnack} because
\begin{align}
h_{\varepsilon}(0,y) & = (67 - \varepsilon) + 2 y + 130 y^2 + 52 y^3 + 32 y^4,
\nonumber \\
h_{\varepsilon}(x,0) & = (67 - \varepsilon) + 2 x + 130 x^2 + 52 x^3 + 32 x^4.
\nonumber 
\end{align}

The polynomial $H$ has exactly four real zeros, which are all positive and given by
\begin{align}
\mathbf{p}_1 & = \left(1 - \sqrt{\frac{1}{8}}, 1 + \sqrt{\frac{1}{8}} \right), 
\; \; \; \; \; \; \;
\mathbf{p}_2 = \left(1 + \sqrt{\frac{1}{8}}, 1 - \sqrt{\frac{1}{8}} \right), \nonumber \\
\mathbf{p}_3 & = \left(\frac{3}{2} - \sqrt{\frac{5}{8}}, \frac{3}{2} + \sqrt{\frac{5}{8}}\right), 
\; \; \; \; \; \; 
\mathbf{p}_4 = \left(\frac{3}{2} + \sqrt{\frac{5}{8}}, \frac{3}{2} - \sqrt{\frac{5}{8}}\right). \label{eq:p_roots}
\end{align}
Since $H \ge 0$ on $\mathbb{R}^2$ and $H(\mathbf{p}_i) = 0$, 
each of the four points from~(\ref{eq:p_roots}) is a global minimum of $H$.
Letting $(x,y) = \mathbf{p}_i + \boldsymbol{\eta}$, 
Taylor-expanding gives $H(\mathbf{p}_i + \boldsymbol{\eta}) 
= (\nabla A(\mathbf{p}_i) \cdot \boldsymbol{\eta})^2
+ 31 (\nabla B(\mathbf{p}_i) \cdot \boldsymbol{\eta})^2 
+ \mathcal{O}(\|\boldsymbol{\eta}\|^3)$. 
The leading-order quadratic form
is positive definite since $\textrm{det}(\nabla A,\nabla B)(\mathbf{p}_i)
= - 8 (p_{i,1} - p_{i,2})(2 p_{i,1} + 2 p_{i,2} - 5) \ne 0$
for every $i = 1,2,3,4$. Hence, all four minima are non-degenerate.

Consider the auxiliary Hamiltonian DS~(\ref{eq:auxiliary_H}).
Since the four points $\mathbf{p}_i$, $i = 1,2,3,4$, 
are non-degenerate minima of $H$, 
they are non-degenerate centers of this DS~\cite{Perko}. 
Therefore, for every sufficiently small $\varepsilon > 0$, 
the level set $\{H = \varepsilon\}$ contains four distinct 
periodic orbits $\Gamma_i$, $i = 1,2,3,4$, each 
surrounding one of the centers and contained in $\mathbb{R}_{>}^2$. 
These periodic orbits are regular ovals of $\{h_{\varepsilon} = 0\}$
because $\{H = \varepsilon\} = \{h_{\varepsilon} = 0\}$.
Furthermore, every point satisfying 
$\nabla H = 2 A \nabla A + 62 B \nabla B = \mathbf{0}$
satisfies $H = 0$, $x = y$ or $x + y = 5/2$. 
Since $H$ is bounded away from zero on both lines, it follows that,
for every sufficiently small $\varepsilon > 0$, 
$\nabla H \ne \mathbf{0}$ on $\{H = \varepsilon\}$,
i.e. $\{h_{\varepsilon} = 0\}$ is regular.
Since a regular quartic has at most $M(4) = 4$ ovals,
$\Gamma_i$, $i = 1,2,3,4$, are the only bounded connected components
of $\{h_{\varepsilon} = 0\}$, and the Harnack bound is attained.
Therefore $\{h_{\varepsilon} = 0\}$ is a chemical Harnack M-curve.
\end{proof}

Lemmas~\ref{lemma:chemical} and~\ref{lemma:cHarnack_4}
jointly give quartic CDS~(\ref{eq:quartic_4})
with four algebraic limit cycles. Its
phase-plane diagram with $\varepsilon = 5/100$ 
is shown in Figure~\ref{fig:1}(b).

\begin{theorem}
\label{theorem:quartic_4}
For every sufficiently small $\varepsilon > 0$, the quartic \emph{CDS}
\begin{align}
\frac{\mathrm{d} x}{\mathrm{d} t} & = 
2 y - 236 x y + 260 y^2 - 340 x^2 y - 680 x y^2 + 156 y^3 
- 368 x^3 y + 2368 x^2 y^2 - 1104 x y^3 + 128 y^4, \nonumber \\
\frac{\mathrm{d} y}{\mathrm{d} t} & = 
(67 - \varepsilon) + 2 x + 130 x^2 - 496 x y + 366 y^2 + 52 x^3 
- 496 x^2 y + 340 x y^2 + 392 y^3 \nonumber \\
& + 32 x^4 - 496 x^3 y + 2288 x^2 y^2 - 2736 x y^3 + 400 y^4,
\label{eq:quartic_4}
\end{align}
admits four hyperbolic algebraic limit cycles 
$\Gamma_1,\Gamma_2,\Gamma_3,\Gamma_4 \subset \{h_{\varepsilon} = 0\} 
\cap \mathbb{R}_{>}^2$, where $h_{\varepsilon}$ is the quartic 
polynomial given in~$(\ref{eq:h_quartic_single})$. 
For each $i = 1,2,3,4$, the cycle $\Gamma_i$ encloses the point $\mathbf{p}_i$,
given in~$(\ref{eq:p_roots})$, and is stable.
Furthermore, $(\ref{eq:quartic_4})$ has no other cycles in $\mathbb{R}_{>}^2$. 
\end{theorem}

\begin{proof}
Let $h_{\varepsilon}$ be the quartic polynomial 
given in~$(\ref{eq:h_quartic_single})$. By Lemma~\ref{lemma:cHarnack_4}, 
$\{h_{\varepsilon} = 0\}$ is a chemical Harnack M-curve 
and hence admits $M(4) = 4$ ovals. Substituting 
$h_{\varepsilon}$ into~(\ref{eq:DS_h_chemical_1}), 
one obtains~(\ref{eq:quartic_4}),
which, by Lemma~\ref{lemma:chemical}, 
is a CDS admitting the stated four hyperbolic limit cycles, 
and no other cycles, in $\mathbb{R}_{>}^2$.
Furthermore, by Lemma~\ref{lemma:cHarnack_4}, each 
$\Gamma_i \subset \{H = \varepsilon\}$ is a sufficiently small
oval enclosing the non-degenerate minimum $\mathbf{p}_i$ 
of $H$ satisfying $H(\mathbf{p}_i) = 0$. Hence,
$h_{\varepsilon} = H - \varepsilon < 0$ inside $\Gamma_i$,
and $h_{\varepsilon} = H - \varepsilon > 0$ immediately outside $\Gamma_i$.
Therefore, by Lemma~\ref{lemma:chemical}, all four limit cycles are stable.
\end{proof}

\noindent 
\textbf{Remark}. A more detailed analysis shows that 
$\{h_{\varepsilon} = 0\}$ constructed in Lemma~\ref{lemma:cHarnack_4}
is a chemical Harnack M-curve for every $0 < \varepsilon < 1/16$. 
Consequently, Theorem~\ref{theorem:quartic_4} 
holds throughout this parameter range.

Theorem~\ref{theorem:quartic_4} gives the following corollary.

\begin{corollary}
\label{corollary:quartic_4}
$H_a^c(4) \ge 4$. 
\end{corollary}

\subsection{Cubic system with two algebraic limit cycles}
\label{sec:cancel_1}
By Lemmas~\ref{lemma:chemical} and~\ref{lemma:cHarnack_3},
there exists a cubic CDS with an algebraic limit cycle.
This attains the maximum number of algebraic limit cycles for cubic 
polynomial DSs with generic invariant 
algebraic curves~\cite{Llibre,Llibre2}.
We now construct a cubic CDS with two algebraic limit cycles. 
Key to this construction is to choose a non-generic 
quartic algebraic curve $\{h = 0\}$ with two ovals 
and appropriate functions $f_1, g_1$ and $c$ such 
that the quartic terms in~(\ref{eq:DS_1}) cancel out. 
A similar approach is used in~\cite{Llibre_Zhao}[Proposition 19]. 

In particular, the authors from~\cite{Llibre_Zhao} 
consider the quartic polynomial 
$\bar{H}(u,v) = -u^2/2 + v^2/2 + u^4/4$ 
whose level set $\{\bar{H} = -1/8\}$ consists purely of two ovals, 
one enclosing $(-1,0)$ and the other enclosing $(1,0)$. 
Let us now define the quartic polynomial 
$H = 8 (\bar{H}+ 1/8)$, given by
\begin{align}
H(u,v) & = 1 - 4 u^2 + 4 v^2 + 2 u^4, \label{eq:H_double}
\end{align}
whose zero-level set $\{H = 0\}$ consists of the same two ovals.
This curve is non-generic in the sense of~\cite{Llibre}
because its highest-degree homogeneous part, $2 u^4$, 
has a repeated factor $u$. 

\begin{lemma}
\label{lemma:cubic_double}
Let $\Gamma_1$ and $\Gamma_2$ be the two ovals of $\{H = 0\}$, 
enclosing $(-1,0)$ and $(1,0)$, respectively, where
$H$ is the quartic polynomial given by~$(\ref{eq:H_double})$.
Consider the cubic \emph{DS}
\begin{align}
\frac{\mathrm{d} u}{\mathrm{d} t} 
& = - 8 c(u,v) v, \nonumber \\
\frac{\mathrm{d} v}{\mathrm{d} t} 
& = 2 \beta u^3 + g_1(u,v) (1 - 4 u^2 + 4 v^2) - 8 c(u,v) u,
\label{eq:cubic_double}
\end{align} 
where, for $\alpha_{0,0}, \alpha_{1,0}, \alpha_{0,1}, \beta \in \mathbb{R}$,  
\begin{align}
g_1(u,v) & = \alpha_{0,0} + \alpha_{1,0} u + \alpha_{0,1} v, \; \; \; \; \; \; \; 
c(u,v) = \frac{1}{4} \left(\beta -  g_1(u,v) u \right).
\label{eq:cubic_double_fgc}
\end{align}
If $c \ne 0$ on $\Gamma_1 \cup \Gamma_2$
and $\alpha_{0,1} \beta \ne 0$, then $\Gamma_1$ and $\Gamma_2$
are hyperbolic algebraic limit cycles of~$(\ref{eq:cubic_double})$.
For each $i = 1,2$, the cycle $\Gamma_i$ is 
stable if $\alpha_{0,1} \beta c < 0$ on $\Gamma_i$, and 
unstable if $\alpha_{0,1} \beta c > 0$ on $\Gamma_i$.
Furthermore, $(\ref{eq:cubic_double})$ 
has no other cycles in $\mathbb{R}^2$.
\end{lemma}

\begin{proof}
Substituting~(\ref{eq:H_double}) into~(\ref{eq:DS_1})
and choosing $f_1 = 0$, one obtains
\begin{align}
\frac{\mathrm{d} u}{\mathrm{d} t} 
& = - 8 c(u,v) v, \nonumber \\
\frac{\mathrm{d} v}{\mathrm{d} t} 
& = 2 u^3 \left(g_1(u,v) u + 4 c(u,v) \right) 
+ g_1(u,v) (1 - 4 u^2 + 4 v^2) - 8 c(u,v) u.
\label{eq:cubic_double_pre}
\end{align} 
The first equation is at most cubic if $c$ is at most quadratic. 
The second equation is then at most cubic if $g_1$ is affine 
and $g_1 u + 4 c$ is constant. 
Therefore we choose~(\ref{eq:cubic_double_fgc}),
under which~(\ref{eq:cubic_double_pre}) coincides 
with~(\ref{eq:cubic_double}). 
Furthermore~(\ref{eq:F}) gives $F(u,v) = \alpha_{0,1} \beta/4$.

Assume that $c \ne 0$ on $\Gamma_1 \cup \Gamma_2$
and $\alpha_{0,1} \beta \ne 0$. For each $i = 1,2$, since $\Gamma_i$ is connected
and $c$ is continuous, $c$ has a constant sign on $\Gamma_i$. 
Given $(u_0,v_0) \in \Omega(\Gamma_i)$, 
the vertical line through this point intersects $\Gamma_i$ at points $(u_0,v_{-})$
and $(u_0,v_{+})$ with $v_{-} \le v_0 \le v_{+}$. Since $c(u_0,\cdot)$ is affine
and hence monotone, $c(u_0,v_0)$ lies between the two boundary values,
which have the same sign. Hence $c$ retains that sign throughout 
$\overline{\Omega(\Gamma_i)}$. 

For each $i = 1,2$, choose a sufficiently small simply connected open set
$\mathbb{S}_i$ containing $\overline{\Omega(\Gamma_i)}$ on which $c \ne 0$
and which is disjoint from the other oval. Since $F \ne 0$ on $\mathbb{R}^2$,
Theorem~\ref{theorem:AL_strong} applies separately in each $\mathbb{S}_i$.
Hence $\Gamma_1$ and $\Gamma_2$ are hyperbolic 
algebraic limit cycles of~(\ref{eq:cubic_double}).
Furthermore, $\{H = 0\}$ consists purely of two disjoint, non-nested ovals, 
so $H \ne 0$ in $\Omega(\Gamma_i)$. Since $H(-1,0) = H(1,0) = -1$, 
$H < 0$ in $\Omega(\Gamma_i)$. By regularity, 
$H$ changes sign across each oval, so $H > 0$ immediately outside it.
Hence $s(\mathbf{x}_{\Gamma_i}) = \textrm{sign}(c(\mathbf{p}_i))$ 
for any $\mathbf{p}_i \in \Gamma_i$. Using~(\ref{eq:K_F}), 
$\textrm{sign}(K_{\Gamma_i}) = \textrm{sign}(\alpha_{0,1} \beta c(\mathbf{p}_i))$,
giving the stated stability. 

Assume for contradiction that there exists an additional cycle 
$\Gamma \not \in \{\Gamma_1,\Gamma_2\}$.
Consider the set $\{c = 0\}$. Since $c(0,v) = \beta/4 \ne 0$, 
the line $u = 0$ does not intersect $\{c = 0\}$.
Since $\alpha_{0,1} \ne 0$, $\{c = 0\} 
= \{(u,v) \in \mathbb{R}^2 \, | \, v = (\beta/u - \alpha_{0,0} 
- \alpha_{1,0} u)/\alpha_{0,1}, \, u \ne 0\}$,
so this set consists of two connected and unbounded branches,
corresponding to $u < 0$ and $u > 0$.
Moreover, $\{H \le 0\} = \overline{\Omega(\Gamma_1)} 
\cup \overline{\Omega(\Gamma_2)}$
does not intersect $\{c = 0\}$, as established above. 
Hence, $H > 0$ on $\{c = 0\}$. Since 
$F = \alpha_{0,1} \beta/4 \ne 0$ is constant, 
it follows from~(\ref{eq:c_derivative}) that
$\mathrm{d} c/\mathrm{d} t \ne 0$ has a constant sign on $\{c = 0\}$. 
Consequently, no cycle intersects $\{c = 0\}$.
Since both branches of $\{c = 0\}$ are connected and unbounded, 
neither can intersect $\Omega(\Gamma)$
without intersecting $\Gamma$. Hence, $c \ne 0$
on $\overline{\Omega(\Gamma)}$. We can therefore choose
a simply connected open set $\mathbb{S}$ containing  
$\overline{\Omega(\Gamma)}$ on which $c \ne 0$. 
Since $F \ne 0$ on $\mathbb{S}$, and $\{H = 0\}$
consists of only two ovals $\Gamma_1$ and $\Gamma_2$,
Theorem~\ref{theorem:AL_strong} excludes any additional
cycle in $\mathbb{S}$, contradicting the existence of $\Gamma$. 
Hence, (\ref{eq:cubic_double}) 
has no other cycles in $\mathbb{R}^2$.
\end{proof}

\noindent 
\textbf{Remark}. \cite{Llibre_Zhao}[Equation~(50)]
is recovered from~(\ref{eq:cubic_double}) by choosing 
$\alpha_{0,0} = \alpha_{1,0} = 0$, $\alpha_{0,1} = 1$
and $\beta = -10$. On $\{H = 0\}$, $2 (u^2 - 1)^2 + 4 v^2 = 1$,
implying $|u| < 2$ and $|v| \le 1/2$, and hence $|u v| < 1$.
Therefore $c(u,v) = -(10 + u v)/4 < 0$, so 
$\alpha_{0,1} \beta c > 0$ on both ovals. 
By Lemma~\ref{lemma:cubic_double},
both limit cycles are unstable,
and no other cycles exist in $\mathbb{R}^2$.

We now construct a cubic CDS with two algebraic limit cycles, 
given by~(\ref{eq:cubic_2}), whose 
phase-plane diagram is shown in Figure~\ref{fig:1}(a).

\begin{theorem}
\label{theorem:cubic_2}
The cubic \emph{CDS}
\begin{align}
\frac{\mathrm{d} x}{\mathrm{d} t} & = 
16 x + 20 y + 58 x^2 -92 x y + 6 y^2 + x^3 - 11 x^2 y + 13 x y^2 + y^3, \nonumber \\
\frac{\mathrm{d} y}{\mathrm{d} t} & = 
2 + 10 x + 17 y + 54 x^2 - 84 x y + 6 y^2 + 5 x^3 - 21 x^2 y + 19 x y^2 + y^3,
\label{eq:cubic_2}
\end{align}
admits two hyperbolic algebraic limit cycles $\Gamma_1, \Gamma_2 
\subset \{h = 0\} \cap \mathbb{R}_{>}^2$, where $h(x,y) = H(x - y - 2, 2 x - 3 y + 2)$ 
and $H$ is the quartic polynomial given by~$(\ref{eq:H_double})$. 
The cycle $\Gamma_1$ encloses $(5,4)$ and is stable, 
while $\Gamma_2$ encloses $(11,8)$ and is unstable.
Furthermore, $(\ref{eq:cubic_2})$ 
has no other cycles in $\mathbb{R}^2$.
\end{theorem}

\begin{proof}
Let $(\alpha_{0,0}, \alpha_{1,0}, \alpha_{0,1}, \beta) = (-8, -3, 1, 1/2)$.
Substitution into~(\ref{eq:cubic_double_fgc}) gives 
\begin{align}
g_1(u,v) & = - 8 - 3 u + v, \; \; \; \; \; \; \; 
c(u,v) = \frac{1}{8} + \frac{1}{4} u (8 + 3 u - v),
\label{eq:cubic_2_chemical}
\end{align}
On the two ovals, $\{H = 0\}$ gives $2 (u^2 - 1)^2 + 4 v^2 = 1$,
so $1/2 < |u| < 3/2$ and $|v| \le 1/2$. Therefore $8 + 3 u - v > 3$. 
Moreover, since $H(0,v) = 1 + 4 v^2 > 0$ for all $v \in \mathbb{R}$, 
the ovals do not intersect the line $u = 0$.
Hence, on the left oval, $u < -1/2$ and $c < 1/8 - 3/8  < 0$,
while on the right oval, $u > 0$ and $c > 0$.
Using also that $\alpha_{0,1} \beta = 1/2 > 0$, 
by Lemma~\ref{lemma:cubic_double},
both ovals are hyperbolic algebraic limit cycles. Furthermore, 
the left limit cycle is stable, while the right one is unstable,
and there are no other cycles. 

Consider the affine change of coordinates
\begin{align}
x & = 8 + 3 u - v, \; \; \; \; y = 6 + 2 u - v.
\end{align}
Since it is invertible, this map preserves the cubic degree, 
cycles and their stability. Furthermore, 
since on both ovals $|u| < 3/2$ and $|v| \le 1/2$,
$x > 3 > 0$ and $y > 5/2 > 0$, i.e. both ovals are located
in the positive quadrant. Moreover, using the inverse map 
$u = -2 + x - y$, $v = 2 + 2 x - 3 y$, DS~(\ref{eq:cubic_double})
with~(\ref{eq:cubic_2_chemical}) maps to~(\ref{eq:cubic_2}),
which is a CDS by Definition~\ref{def:CDS}.
\end{proof}

Theorem~\ref{theorem:cubic_2} gives the following corollary.

\begin{corollary}
\label{corollary:cubic_2}
$H_a^c(3) \ge 2$. 
\end{corollary}

\section{Algebraic limit cycles in reversible chemical dynamical systems}
\label{sec:5}

In Section~\ref{sec:4}, we have studied (semi-)chemical polynomials 
and used the results from Section~\ref{sec:3}
to construct general CDSs with prescribed algebraic limit cycles. 
In this section, we specialize this approach to 
reversible CDSs~\cite{Feinbergbook,Janos,YuCraciun_2018}, 
which are subject to additional constraints on their coefficients. 

\begin{definition}
$($\textbf{Reversible \emph{CDS}}$)$ 
\label{def:reversibility}
Let~$(\ref{eq:DS})$ be a \emph{CDS}. A \emph{CRN} realization $\mathcal{R}$
of~$(\ref{eq:DS})$ is said to be \emph{reversible} if, for every reaction 
$i X + j Y \xrightarrow[]{\alpha} i' X + j' Y$ in $\mathcal{R}$,
the reverse reaction 
$i' X + j' Y  \xrightarrow[]{\alpha'} i X + j Y$ is also in $\mathcal{R}$. 
We say that \emph{CDS}~$(\ref{eq:DS})$ is \emph{reversible}
if it admits a reversible \emph{CRN} realization.  
\end{definition}

\noindent
\textbf{Remark}. The notion of reversible CDSs introduced in 
Definition~\ref{def:reversibility} is distinct from 
the notion of reversible dynamical systems, which refers to
time-reversal symmetry. 
 
For any integer $n \ge 1$, let us define the triangular lattice
\begin{align}
T_n & = \{(i,j) \in \mathbb{Z}_{\ge}^2 \, | \, i + j \le n\}.
\label{eq:buffering_triangle}
\end{align}
Assume that a CRN $\mathcal{R}$ contains only reactions of the form 
$i X + j Y \xrightarrow[]{\alpha} i' X + j' Y$, where 
$i + j \le n$ and $i' + j' \le n$. Then this CRN
can be represented as a directed graph on $T_n$
called a \emph{Euclidean embedded graph} (E-graph)~\cite{YuCraciun_2018}.
In particular, a chemical reaction 
$i X + j Y \xrightarrow[]{\alpha} i' X + j' Y$ in $\mathcal{R}$
is represented as a directed edge (with weight $\alpha$)
from the vertex $(i,j) \in T_n$ 
to the vertex $(i',j') \in T_n$. 
In this representation, a reversible CRN is a graph on $T_n$ 
in which every directed edge has its reverse. 

\begin{figure}[!htbp]
\vskip 0.0cm
\leftline{\hskip 0.3cm (a) E-graph of~(\ref{eq:example_reversible_1})
\hskip  2.0cm (b) E-graph of~(\ref{eq:Escher})
\hskip  2.0cm (c) E-graph of~(\ref{eq:cubic_R})}
\vskip 0.0cm
\centerline{
\hskip -1.0cm
\includegraphics[width=0.3\columnwidth]{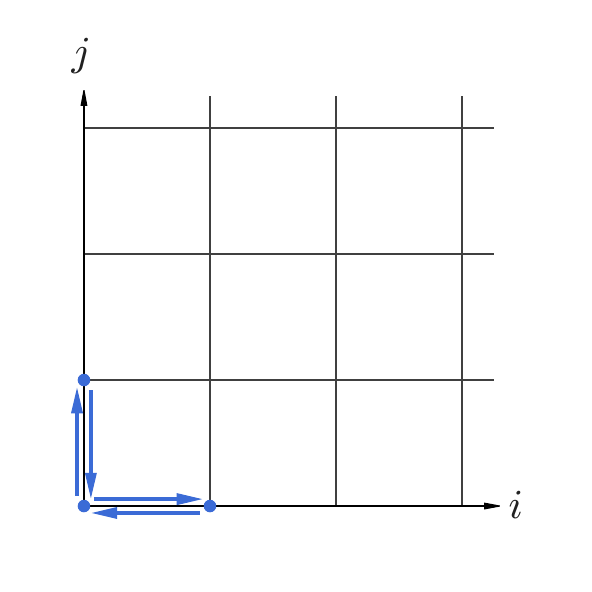}
\hskip 0.3cm
\includegraphics[width=0.3\columnwidth]{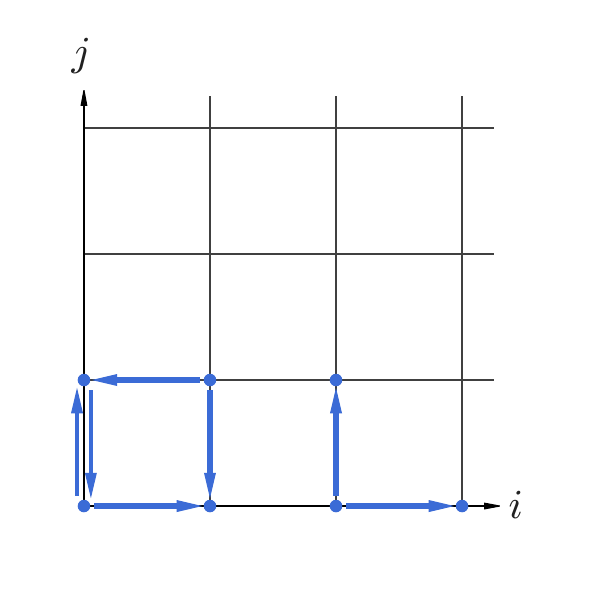}
\hskip 0.3cm
\includegraphics[width=0.3\columnwidth]{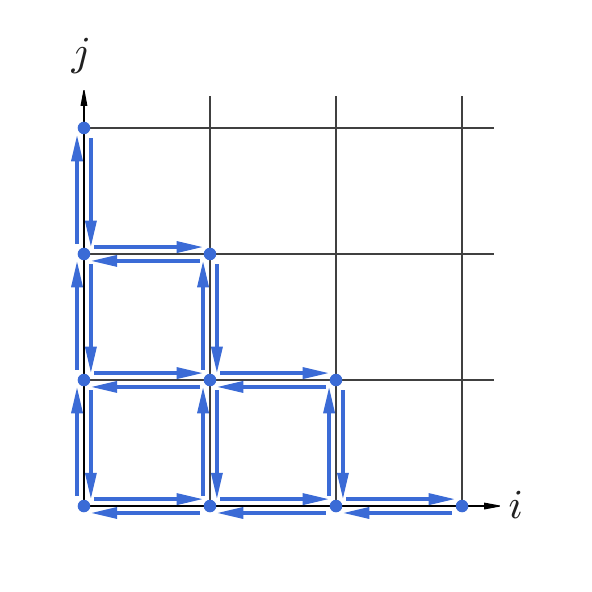}
}
\vskip 0.0cm
\caption{\it Panels \emph{(a)}--\emph{(c)} 
display E-graphs compatible with 
\emph{CDSs}~$(\ref{eq:example_reversible_1})$,
$(\ref{eq:Escher})$ and~$(\ref{eq:cubic_R})$, 
respectively.} 
\label{fig:2}
\end{figure}

\begin{example}
The \emph{CDS}
\begin{align}
\frac{\mathrm{d} x}{\mathrm{d} t} & = 1 - x, \nonumber \\
\frac{\mathrm{d} y}{\mathrm{d} t} & = 1 - y, 
\label{eq:example_reversible_1}
\end{align}
is reversible because it admits a reversible \emph{CRN} 
\begin{align}
\varnothing &\xrightleftharpoons[1]{1} X, \; \; \; 
\varnothing \xrightleftharpoons[1]{1} Y.
\end{align}
Its (non-weighted) E-graph on $T_1 = \{(0,0),(1,0),(0,1)\}$ 
is shown in \emph{Figure~\ref{fig:2}(a)}. 
\end{example}

\begin{definition}
\label{def:cross_reversibility}
$($\textbf{Cross-stencil reversibility}$)$ 
A \emph{CDS} of degree $n$ is \emph{cross-stencil reversible}
if it admits a reversible \emph{CRN} realization whose E-graph on $T_n$ 
is such that every vertex is connected to all its
horizontal and vertical nearest neighbors within $T_n$
by edges with positive weights, and to no other vertices. 
\end{definition}

CDS~(\ref{eq:example_reversible_1}) is cross-stencil 
reversible on $T_1$, while Figure~\ref{fig:2}(c)
shows a cross-stencil reversible CDS on $T_3$.
The edges form a cross at each interior vertex,
which is truncated at the boundary.

We now present conditions, taking form of linear inequalities 
on the polynomial coefficients $\alpha_{i,j}$ and $\beta_{i,j}$,
that are necessary and sufficient for CDS~(\ref{eq:DS}) to belong to this subclass.

\begin{lemma}
\label{lemma:reversible_1}
The \emph{DS}~$(\ref{eq:DS})$ of degree $n \ge 1$
is a cross-stencil reversible \emph{CDS}
if and only if its coefficients satisfy the following linear inequalities:
\begin{align}
\alpha_{0,j} & > 0, \; \; \; j = 0,1,\ldots,n-1,
\; \; \text{ and } \; \; \alpha_{0,n} = 0, \nonumber \\
\beta_{i,0} & > 0, \; \; \; i = 0,1,\ldots,n-1,
\; \; \text{ and } \; \; \beta_{n,0} = 0, \nonumber \\
\alpha_{i,n-i}& < 0, \; \; \; i = 1,2,\ldots,n, \nonumber \\
\beta_{i,n-i} & < 0, \; \; \; i = 0,1,\ldots,n-1.
\label{eq:reversible_1}
\end{align}
\end{lemma}

\begin{proof}
Consider the reversible CRN consisting of the horizontal reactions
\begin{align}
i X + j Y & \xrightleftharpoons[h_{i,j}^{-}]{h_{i,j}^{+}} (i + 1) X + j Y, 
\; \; \; 0 \le i + j \le n - 1,
\end{align}
and the vertical reactions 
\begin{align}
i X + j Y & \xrightleftharpoons[v_{i,j}^{-}]{v_{i,j}^{+}} i X + (j + 1) Y, 
\; \; \; 0 \le i + j \le n - 1,
\end{align}
where all the rate coefficients $h_{i,j}^{\pm}, v_{i,j}^{\pm}$ are positive.
The corresponding CDS is given by
\begin{align}
\frac{\mathrm{d} x}{\mathrm{d} t} & = \sum_{i + j = 0}^{n-1} 
\left(h_{i,j}^{+} x^i y^j - h_{i,j}^{-} x^{i+1} y^j \right), \nonumber \\
\frac{\mathrm{d} y}{\mathrm{d} t} & = \sum_{i + j = 0}^{n-1} 
\left(v_{i,j}^{+} x^i y^j - v_{i,j}^{-} x^i y^{j+1} \right).
\label{eq:comparison_1}
\end{align}
By Definition~\ref{def:cross_reversibility}, this CDS 
is cross-stencil reversible.

Assume first that DS~(\ref{eq:DS}) satisfies conditions~(\ref{eq:reversible_1}).
We now show that the coefficients $h_{i,j}^{\pm}, v_{i,j}^{\pm}$
can be chosen so that~(\ref{eq:comparison_1}) and~(\ref{eq:DS}) coincide. 
Let us rewrite the right-hand side of 
the first equation in~(\ref{eq:comparison_1}) as
\begin{align}
\sum_{i + j = 0}^{n-1} \left(h_{i,j}^{+} x^i y^j - h_{i,j}^{-} x^{i+1} y^j \right)
& = \sum_{i + j = 0}^{n-1} h_{i,j}^{+} x^i y^j 
- \sum_{i + j = 0}^{n-1} h_{i,j}^{-} x^{i+1} y^j
= \sum_{i + j = 0}^{n-1} h_{i,j}^{+} x^i y^j 
- \sum_{\substack{i + j = 1, \\ i \ge 1}}^{n} h_{i-1,j}^{-} x^{i} y^j \nonumber \\
& = \sum_{j = 0}^{n-1} h_{0,j}^{+} y^j 
+ \sum_{\substack{i + j = 1, \\ i \ge 1}}^{n-1} 
\left(h_{i,j}^{+} - h_{i-1,j}^{-} \right) x^{i} y^j
- \sum_{i = 1}^n h_{i-1,n-i}^{-} x^{i} y^{n-i}.
\end{align}
Equating this with $\sum_{i + j = 0}^n \alpha_{i,j} x^i y^j$, one obtains
\begin{align}
h_{0,j}^{+} & = \alpha_{0,j}, 
\; \; \; \; \; \; \; \; \; \; \; \;  j = 0,1,\ldots,n-1, \nonumber \\
h_{i,j}^{+} - h_{i-1,j}^{-} & = \alpha_{i,j}, 
\; \; \; \; \; \; \; \; \; \; \; \; i \ge 1, \; \; i + j \le n - 1, \nonumber \\
h_{i-1,n-i}^{-} & = - \alpha_{i,n-i}, 
\; \; \; \; \; \;  i = 1,2,\ldots, n, \nonumber \\
\alpha_{0,n} & = 0. \label{eq:comparison_equations}
\end{align}
By assumption~(\ref{eq:reversible_1}), 
one can choose $h_{0,j}^{+} = \alpha_{0,j} > 0$ for $j = 0,1,\ldots,n-1$ and
$h_{i-1,n-i}^{-} = - \alpha_{i,n-i} > 0$ for $i = 1,2,\ldots, n$.
For the remaining coefficients,
$i \ge 1, i + j \le n - 1$, 
choose $h_{i-1,j}^{-} > \max\{0,-\alpha_{i,j}\} \ge 0$
and $h_{i,j}^{+} = \alpha_{i,j} + h_{i-1,j}^{-} > 0$. 
Hence, the first equation in~(\ref{eq:DS})
can be reproduced with the first equation in~(\ref{eq:comparison_1}) .
An analogous argument, equating the second equations 
in~(\ref{eq:comparison_1}) and~(\ref{eq:DS}), proves
the same for the second equations. Therefore~(\ref{eq:DS})
is cross-stencil reversible.

Conversely, assume that~(\ref{eq:DS})
is cross-stencil reversible.
Then it admits a CRN realization whose corresponding
CDS is given by~(\ref{eq:comparison_1}), with 
all the rate coefficients $h_{i,j}^{\pm}, v_{i,j}^{\pm}$ positive. 
Comparing the first equations in~(\ref{eq:comparison_1}) and~(\ref{eq:DS})
yields~(\ref{eq:comparison_equations}), which imply 
the conditions on $\alpha_{i,j}$ in~(\ref{eq:reversible_1}). 
Analogous argument yields the conditions on $\beta_{i,j}$. 
\end{proof}

\begin{example}
\emph{CDS}~$(\ref{eq:Escher})$ is not 
cross-stencil reversible
since $\beta_{2,0} = 5/2 \ne 0$, 
violating~$(\ref{eq:reversible_1})$.
The E-graph of its 
\emph{CRN} realization~$(\ref{eq:Escher_CRN})$
is shown in \emph{Figure~\ref{fig:2}(b)}. 
\end{example}

\subsection{Dynamics induced by buffered polynomials}
To systematically construct reversible CDSs with algebraic limit cycle, 
we introduce the notion of buffered polynomials 
that have positive coefficients on the boundary of 
a triangular lattice, while the coefficients 
in the interior remain unrestricted. 

\begin{definition}
$($\textbf{Buffered polynomial}$)$ 
\label{def:buffered}
The $n$-degree polynomial~$(\ref{eq:polynomial_h})$ is said to be 
\emph{buffered} on $T_n$ if $\kappa_{i,j} > 0$
for every $(i,j) \in \partial T_n$, where $\partial T_n$
denotes the boundary lattice points of~$(\ref{eq:buffering_triangle})$. 
Equivalently,
\begin{align}
\kappa_{i,0} & > 0 \text{ for every } i = 0,1,\ldots,n, \nonumber \\
\kappa_{0,j} & > 0 \text{ for every } j = 1,2,\ldots,n, \nonumber \\
\kappa_{i,n - i} & > 0 \text{ for every } i = 1,2,\ldots,n-1.
\label{eq:buffered}
\end{align}
\end{definition}

\noindent
\textbf{Remark}. By Definitions~\ref{def:semi_chemical} and~\ref{def:buffered}, 
every buffered polynomial is chemical, but the converse is false,
i.e. buffering is a stronger condition than chemicality.

We now use Lemma~\ref{lemma:reversible_1} and Theorem~\ref{theorem:AL_strong}
to prove the following result. 

\begin{lemma}
\label{lemma:CDS_R_1}
Let $h : \mathbb{R}^2 \to \mathbb{R}$
be a buffered polynomial of degree $n \ge 1$
such that $\{h = 0\}$ contains $N$ 
regular simple closed curves
$\Gamma_1, \Gamma_2, \ldots, 
\Gamma_N  \subset \mathbb{R}_{>}^2$. 
Then, for every sufficiently small $\varepsilon > 0$, 
\begin{align}
\frac{\mathrm{d} x}{\mathrm{d} t} 
& = (1 - x) h(x,y) 
- \varepsilon \frac{\partial h(x,y)}{\partial y}, \nonumber \\
\frac{\mathrm{d} y}{\mathrm{d} t} 
& = (1 - y) h(x,y) 
+ \varepsilon \frac{\partial h(x,y)}{\partial x},
\label{eq:CDS_R_1}
\end{align}
is a cross-stencil reversible
\emph{CDS} of degree $n+1$.
Furthermore, $(\ref{eq:CDS_R_1})$ has $N$ hyperbolic 
algebraic limit cycles 
$\Gamma_1, \Gamma_2, \ldots, \Gamma_N$.
For each $i = 1,2,\ldots, N$, the cycle $\Gamma_i$
is stable (respectively, unstable) if 
$h(\mathbf{p}_i) > 0$ (respectively, $h(\mathbf{p}_i) < 0$),
where $\mathbf{p}_i$ is any point
sufficiently close to $\Gamma_i$ and in 
the exterior of $\Gamma_i$.
Furthermore, $(\ref{eq:CDS_R_1})$
has no other cycles in $\mathbb{R}_{>}^2$
if $\{h = 0\}$ contains no 
bounded connected components in $\mathbb{R}_{>}^2$
other than $\Gamma_1, \Gamma_2, \ldots, \Gamma_N$.
\end{lemma}

\begin{proof}
Writing $h$ as in~(\ref{eq:polynomial_h}),
and substituting into~(\ref{eq:CDS_R_1}), one obtains
\begin{align}
\frac{\mathrm{d} x}{\mathrm{d} t} 
& = \sum_{i+j = 0}^n \kappa_{i,j} x^i y^j
- \sum_{\substack{i + j = 1, \\ i \ge 1}}^{n+1} \kappa_{i - 1,j} x^i y^j
- \varepsilon \sum_{i+j = 0}^{n-1} (j+1) \kappa_{i,j+1} x^i y^j, \nonumber \\
\frac{\mathrm{d} y}{\mathrm{d} t} 
& = \sum_{i+j = 0}^n \kappa_{i,j} x^i y^j
- \sum_{\substack{i + j = 1, \\ j \ge 1}}^{n+1} \kappa_{i,j-1} x^{i} y^{j}
+ \varepsilon \sum_{i+j = 0}^{n-1} (i + 1) \kappa_{i+1,j} x^{i} y^{j}.
\label{eq:CDS_R_1_expanded}
\end{align}
Writing the two right-hand sides as $\sum_{i + j = 0}^{n+1} \alpha_{i,j} x^i y^j$
 and $\sum_{i + j = 0}^{n+1} \beta_{i,j} x^i y^j$, 
respectively, it follows that
\begin{align}
\alpha_{0,j} & = \kappa_{0,j} - (j + 1) \varepsilon \kappa_{0,j+1}, 
\; \; \; j = 0,1,2,\ldots, n - 1, \; \; \; \alpha_{0,n} = \kappa_{0,n}, 
\; \; \; \alpha_{0,n+1} = 0, \nonumber \\
\beta_{i,0} & = \kappa_{i,0} + (i + 1) \varepsilon \kappa_{i+1,0}, 
\; \; \; \; i = 0,1,2,\ldots, n - 1, \; \; \; \; \beta_{n,0} = \kappa_{n,0}, 
\; \; \; \beta_{n+1,0} = 0, \nonumber \\
\alpha_{i,n+1-i} & = - \kappa_{i-1,n+1-i}, \; \; \; i = 1,2,\ldots, n + 1, \nonumber \\
\beta_{i,n+1-i} & = - \kappa_{i,n-i}, 
\; \; \; \; \; \; \; \; \; \;  i = 0,1,\ldots,n. 
\label{eq:CDS_R_1_expanded_coeff}
\end{align}
Define 
\begin{align}
\varepsilon_0 & = \min_{0 \le j \le n - 1} 
\frac{\kappa_{0,j}}{(j+1) \kappa_{0,j+1}}. \label{eq:eps_1}
\end{align}
Since $h$ is buffered, (\ref{eq:buffered}) holds
and hence $\varepsilon_0 > 0$.
Furthermore, for every $0 < \varepsilon < \varepsilon_0$, 
the coefficients in~(\ref{eq:CDS_R_1_expanded_coeff})
satisfy the conditions~(\ref{eq:reversible_1}).
Therefore, by Lemma~\ref{lemma:reversible_1},
(\ref{eq:CDS_R_1}) is a cross-stencil reversible CDS.

CDS~(\ref{eq:CDS_R_1}) is of the form~(\ref{eq:DS_1}) with 
$f_1(x,y) = (1 - x)$, $g_1(x,y) = (1 - y)$ and $c(x,y) = \varepsilon$.
Since, for every $\varepsilon > 0$, 
$c(x,y) > 0$, $F(x,y) = - 2 \varepsilon < 0$
and $\textrm{sign}(K_{\Gamma_i}) = -\textrm{sign}(h(\mathbf{p}_i))$,
where $\mathbf{p}_i$ is any point immediately outside $\Gamma_i$,
Theorem~\ref{theorem:AL_strong} with 
$\mathbb{S} = \mathbb{R}_{>}^2$ yields the remaining conclusions.
\end{proof}

\noindent
\textbf{Remark}. The conclusions of Lemma~\ref{lemma:CDS_R_1}
hold for every $0 < \varepsilon < \varepsilon_0$, where
$\varepsilon_0$ is given by~(\ref{eq:eps_1}).

\subsubsection{Construction of buffered polynomials} 
The following is a natural generalization of Definition~\ref{def:chemical_Harnack}.

\begin{definition}
\label{def:buffered_Harnack}
$($\textbf{Buffered Harnack M-curves}$)$ 
Let $\{h = 0\}$ be a Harnack M-curve of degree $n$.
If every bounded connected component of $\{h = 0\}$
is contained in the positive quadrant $\mathbb{R}_{>}^2$
and $h$ is buffered on $T_n$, 
then $\{h = 0\}$ is called a \emph{buffered Harnack M-curve}.
\end{definition} 

\begin{lemma}
\label{lemma:bHarnack_2}
There exists no buffered Harnack M-curve of degree $n = 2$. 
\end{lemma}

\begin{proof}
The conclusion follows directly from Lemma~\ref{lemma:cHarnack_2}.
\end{proof}

\begin{lemma}
\label{lemma:bHarnack_3}
There exists a buffered Harnack M-curve of degree $n = 3$. 
\end{lemma}

\begin{proof}
Let $h_{\varepsilon}$ be the polynomial from 
Lemma~\ref{lemma:cHarnack_3} and fix
$0 < \varepsilon < 1$ so that $\{h_{\varepsilon} = 0\}$
is a chemical Harnack M-curve of degree $3$.
Define the auxiliary cubic polynomial 
\begin{align}
p(x,y) & = x + x^2 + y + y^2 + x y^2 + x^2 y,
\end{align}
and, for every $\delta > 0$, consider
\begin{align}
q_{\delta}(x,y) = h_{\varepsilon}(x,y) + \delta p(x,y)
= (1 - \varepsilon) + \delta x + \delta y + \delta x^2 - 3 x y + \delta y^2 
+ x^3 + \delta x y^2 + \delta x^2 y +  y^3. \nonumber 
\end{align}
For every $\delta > 0$, $q_{\delta}$
is buffered on $T_3$ by Definition~\ref{def:buffered}.
Furthermore, since Harnack M-curves of degree $n$ 
persist under sufficiently small polynomial perturbations
of degree at most $n$, it follows that
$\{q_{\delta} = 0\}$ is a cubic Harnack M-curve
for every sufficiently small $\delta > 0$.
\end{proof}

The quartic chemical Harnack M-curve constructed in Lemma~\ref{lemma:cHarnack_4}
contains monomials $x^3 y$ and $x y^3$ with negative coefficients. 
Consequently, the perturbation approach used in Lemma~\ref{lemma:bHarnack_3}
cannot be applied to obtain a buffered counterpart of the same degree.
We leave the existence and construction of buffered Harnack M-curves of degree $n \ge 4$ 
as an open problem (see also Section~\ref{sec:6}). 
Instead, we now present a systematic method for constructing 
general buffered polynomials. 

By Lemmas~\ref{lemma:semichemical_construct} and~\ref{lemma:chemical_construct}, 
every polynomial of degree $n$ can be transformed
into a semi-chemical polynomial of degree $n$
and a chemical polynomial of degree $n+1$
with all its regular ovals preserved and contained in $\mathbb{R}_{>}^2$, 
respectively. We now show that a buffered polynomial of degree $n+2$
can be constructed with the same properties.

\begin{lemma}
\label{lemma:buffered_construct}
Let $H : \mathbb{R}^2 \to \mathbb{R}$ 
be a polynomial of degree $n$ such that
$\{H = 0\}$ contains $N \ge 1$ regular simple closed curves.
Then there exists a buffered
polynomial of degree $n + 2$ 
whose zero set contains $N$ regular 
simple closed curves in $\mathbb{R}_{>}^2$.
\end{lemma}

\begin{proof}
Let $H$ be a polynomial of degree $n$ such that
$\{H = 0\}$ contains $N \ge 1$ regular simple closed curves 
$\Gamma_1, \Gamma_2, \ldots, \Gamma_N$.
Let $H_n$ denote the $n$-degree homogeneous part of $H$. Since $H_n \not \equiv 0$, 
there exists a vector $\mathbf{w} \in \mathbb{R}^2$ such that $H_n(\mathbf{w}) \ne 0$. 
Replacing $H$ by $-H$ if necessary, we may assume that $H_n(\mathbf{w}) > 0$.

Consider the homogeneous polynomial $H_n(x \mathbf{w}_1 + y \mathbf{w}_2)$, 
where $\mathbf{w}_1, \mathbf{w}_2 \in \mathbb{R}^2$.  
If $\mathbf{w}_1 = \mathbf{w}_2 = \mathbf{w}$, then, 
by homogeneity and the binomial theorem, 
\begin{align}
H_n(x \mathbf{w}_1 + y \mathbf{w}_2) 
& = H_n((x + y) \mathbf{w})
= H_n(\mathbf{w}) (x + y)^n
= H_n(\mathbf{w}) \sum_{i = 0}^n \binom{n}{i} x^i y^{n-i}.
\end{align}
Hence, all $n+1$ coefficients of this polynomial are positive. 
Since the coefficients of $H_n(x \mathbf{w}_1 + y \mathbf{w}_2)$
depend continuously on $\mathbf{w}_1$ and $\mathbf{w}_2$, they
remain positive for every $\mathbf{w}_1$ and $\mathbf{w}_2$
sufficiently close to $\mathbf{w}$. 

Choose any two such vectors $\mathbf{w}_1$ and $\mathbf{w}_2$
that are linearly independent, and define the invertible
linear map $L(x,y) = x \mathbf{w}_1 + y \mathbf{w}_2$. 
Then the polynomial defined by $\tilde{h}(x,y) = H(L(x,y))$
is of degree $n$, and all the $n+1$ coefficients of its
$n$-degree homogeneous part are positive. Furthermore, the curves
$\tilde{\Gamma}_i = L^{-1}(\Gamma_i)$, $i = 1,2,\ldots,N$, 
are regular simple closed curves contained in $\{\tilde{h} = 0\}$. 

By Lemma~\ref{lemma:h_construct}, 
for every sufficiently large $a, b > 0$,
the $(n+2)$-degree polynomial $\hat{h}(x,y) = x y \tilde{h}(x - a, y - b)$
is chemical, and $(\tilde{\Gamma}_i + (a,b))$, 
$i = 1,2,\ldots,N$, are regular simple closed curves in
$\{\hat{h} = 0\}$ contained in $\mathbb{R}_{>}^2$. 
Writing $\hat{h}(x,y) = \sum_{i+j = 0}^{n+2} \gamma_{i,j} x^i y^j$,
we have, by construction, $\gamma_{i,j} = 0$
if $i = 0$ or $j = 0$, and $\gamma_{i,n+2-i} > 0$
for every $i = 1,2,\ldots,n+1$. Geometrically, the coefficients of $\hat{h}$
on the diagonal side of $T_{n+2}$, excluding its endpoints, are positive, 
while the coefficients on the coordinate axes
of $T_{n+2}$ are zero. 

Fix such $a,b > 0$, and let
\begin{align}
p(x,y) & = 1 + \sum_{i = 1}^{n+2} x^i + \sum_{j = 1}^{n+2} y^j.
\end{align}
Define the perturbed polynomial of degree $n + 2$ by
\begin{align}
h_{\varepsilon}(x,y) & = \hat{h}(x,y) + \varepsilon p(x,y).
\end{align}
By the above properties of $\hat{h}$, for every $\varepsilon > 0$,
all the coefficients of $h_{\varepsilon}(x,y)$
on $\partial T_{n+2}$ are positive, and hence $h_{\varepsilon}$ 
is buffered on $\partial T_{n+2}$ by Definition~\ref{def:buffered}.
Furthermore, since the curves $(\tilde{\Gamma}_i + (a,b))$, $i = 1,2,\ldots,N$,
are regular and pairwise disjoint, the implicit function theorem and compactness
imply that, for every sufficiently small $\varepsilon > 0$, 
$\{h_{\varepsilon} = 0\}$ contains $N$ pairwise disjoint
regular simple closed curves 
$\Gamma_{i,\varepsilon} \subset \mathbb{R}_{>}^2$ satisfying 
$d_H(\Gamma_{i,\varepsilon},\tilde{\Gamma}_i + (a,b)) \to 0$ 
as $\varepsilon \to 0$
for every $i = 1,2,\ldots,N$, completing the proof.
\end{proof}

\begin{theorem}
\label{theorem:buffered_n}
For every integer $n \ge 4$, $H_{a,r}^c(n) \ge M(n-3)$,
where $M(n)$ is given by~$(\ref{eq:M_n})$. 
\end{theorem}

\begin{proof}
For $n = 4$, the inequality follows from $M(1) = 0$. Assume that $n \ge 5$.
Let $H$ be a polynomial of degree $n - 3$ such that $\{H = 0\}$
is a Harnack M-curve with $M(n - 3)$ ovals. 
By Lemma~\ref{lemma:buffered_construct}, there exists a buffered 
polynomial $h_{\varepsilon}$ of degree $n - 1$ such that 
$\{h_{\varepsilon} = 0\}$ contains $M(n - 3)$ regular ovals 
in $\mathbb{R}_{>}^2$. Therefore, by Lemma~\ref{lemma:CDS_R_1},
the cross-stencil reversible CDS~(\ref{eq:CDS_R_1}) of degree $n$ admits $M(n - 3)$
hyperbolic algebraic limit cycles in $\mathbb{R}_{>}^2$,
giving $H_{a,r}^c(n) \ge M(n-3)$.
\end{proof}

\subsection{Cubic system with an algebraic limit cycle}
Since buffered Harnack M-curves of degree two do not exist, 
Lemma~\ref{lemma:CDS_R_1} does not ensure existence of
a cubic reversible CDS with an algebraic limit cycle. 
To establish such existence, we add suitable cubic perturbations 
to~(\ref{eq:Escher}) and invoke Theorem~\ref{theorem:perturbations}
to obtain the family of reversible cubic CDSs~(\ref{eq:cubic_R}),
each possessing identical elliptical limit cycle. 
All members of this family are 
cross-stencil reversible, with the common 
E-graph shown in Figure~\ref{fig:2}(c).
The phase-plane diagram of the member with $\mu = 1/100$ 
is shown in Figure~\ref{fig:1}(c).

\begin{theorem}
\label{theorem:cubic_R}
For every sufficiently small $\mu > 0$, 
\begin{align}
\frac{\mathrm{d} x}{\mathrm{d} t} & = 
\left(\frac{3}{2} - 32 \mu \right) + \frac{49}{2} \mu x
+ \frac{7}{2} \mu y + (2 + 54 \mu) x^2 
- (1 + 26 \mu) x y + 4 \mu y^2
- \mu x^3 - \mu x^2 y - \mu x y^2, \nonumber \\
\frac{\mathrm{d} y}{\mathrm{d} t} & = 
\left(\frac{17}{4} - 40 \mu \right) + 4 \mu x 
- \left(1 + \frac{7}{2} \mu \right) y + \left(\frac{5}{2} + 80 \mu \right) x^2
- (1 + 38 \mu) x y + 13 \mu y^2
- \mu x^2 y - \mu x y^2 - \mu y^3,
\label{eq:cubic_R}
\end{align}
is a cross-stencil reversible
\emph{CDS} of degree $3$. Furthermore, 
$\Gamma = \{h = 0\} \subset \mathbb{R}_{>}^2$ 
is a stable hyperbolic algebraic limit cycle,
where $h(x, y) = 19 + 20 x - 16 y + 10 x^2 + 4 y^2 - 12 x y$.
\end{theorem}

\begin{proof}
CDS~(\ref{eq:Escher}) can be written in the form~(\ref{eq:DS_1}) with 
$f_1,g_1,c$ and $h$ given in~(\ref{eq:Escher_functions}).
Consider the perturbed cubic system
\begin{align}
\frac{\mathrm{d} x}{\mathrm{d} t} 
& = \left(f_1(x,y) + \mu f_2(x,y) \right) h(x,y) - 
\left(c(x,y) + \mu d(x,y) \right) \frac{\partial h(x,y)}{\partial y}, 
\nonumber \\
\frac{\mathrm{d} y}{\mathrm{d} t} 
& = \left(g_1(x,y) + \mu g_2(x,y) \right) h(x,y) 
+ \left(c(x,y) + \mu d(x,y) \right) \frac{\partial h(x,y)}{\partial x},
\label{eq:Escher_perturbed}
\end{align} 
where $\mu > 0$, and 
\begin{align}
f_2(x,y) & = \bar{\alpha}_{0,0} + \bar{\alpha}_{1,0} x + \bar{\alpha}_{0,1} y, \; \; \; 
g_2(x,y) = \bar{\beta}_{0,0} + \bar{\beta}_{1,0} x + \bar{\beta}_{0,1} y, \nonumber \\ 
d(x,y) & = \bar{\gamma}_{0,0} + \bar{\gamma}_{1,0} x + \bar{\gamma}_{0,1} y
+ \bar{\gamma}_{2,0} x^2 + \bar{\gamma}_{1,1} x y + \bar{\gamma}_{0,2} y^2.
\label{eq:cubic_R_fgd}
\end{align}
Writing the right-hand sides of~(\ref{eq:Escher_perturbed}) 
as $\sum_{i + j = 0}^{3} \alpha_{i,j} x^i y^j$
and $\sum_{i + j = 0}^{3} \beta_{i,j} x^i y^j$, respectively, 
cross-stencil reversibility conditions~(\ref{eq:reversible_1}) yield, 
for every fixed $\mu > 0$, a system of $14$ linear constraints
in the $12$ coefficients of $f_2,g_2,d$, given by
\begin{align}
\frac{3}{2} + \mu (19 \bar{\alpha}_{0,0} + 16 \bar{\gamma}_{0,0}) > 0, \; \; \; 
-16 \bar{\alpha}_{0,0} + 19 \bar{\alpha}_{0,1} 
- 8 \bar{\gamma}_{0,0} + 16 \bar{\gamma}_{0,1} > 0, \nonumber \\
\bar{\alpha}_{0,0} - 4 \bar{\alpha}_{0,1} - 2 \bar{\gamma}_{0,1} 
+ 4 \bar{\gamma}_{0,2} > 0, \; \; \; 
\bar{\alpha}_{0,1} - 2 \bar{\gamma}_{0,2} = 0, \nonumber \\
\frac{17}{4} + \mu (19 \bar{\beta}_{0,0} + 20 \bar{\gamma}_{0,0}) > 0, \; \; \; 
20 \bar{\beta}_{0,0} + 19 \bar{\beta}_{1,0} 
+ 20 \bar{\gamma}_{0,0} + 20 \bar{\gamma}_{1,0} > 0, \nonumber \\
\frac{5}{2} + 10 \mu (\bar{\beta}_{0,0} + 2 \bar{\beta}_{1,0}
+ 2 \bar{\gamma}_{1,0} + 2 \bar{\gamma}_{2,0}) > 0, \; \; \; 
\bar{\beta}_{1,0} + 2 \bar{\gamma}_{2,0} = 0, \nonumber \\
- 3 \bar{\alpha}_{0,1} + \bar{\alpha}_{1,0} 
+ 3 \bar{\gamma}_{0,2} - 2 \bar{\gamma}_{1,1} < 0, \; \; \; 
5 \bar{\alpha}_{0,1} - 6 \bar{\alpha}_{1,0}
+ 6 \bar{\gamma}_{1,1} - 4 \bar{\gamma}_{2,0} < 0, \; \; \; 
5 \bar{\alpha}_{1,0} + 6 \bar{\gamma}_{2,0} < 0,\nonumber \\
\bar{\beta}_{0,1} - 3 \bar{\gamma}_{0,2}  < 0, \; \; \; 
- 3 \bar{\beta}_{0,1} + \bar{\beta}_{1,0}
+ 5 \bar{\gamma}_{0,2} - 3 \bar{\gamma}_{1,1} < 0, \; \; \; 
5 \bar{\beta}_{0,1} - 6 \bar{\beta}_{1,0}
+ 10 \bar{\gamma}_{1,1} - 6 \bar{\gamma}_{2,0} < 0.
\label{eq:R_1_cubic}
\end{align}
The two equalities give $\bar{\gamma}_{0,2} = \bar{\alpha}_{0,1}/2$
and $\bar{\gamma}_{2,0} = - \bar{\beta}_{1,0}/2$.
Substituting these identities into~(\ref{eq:R_1_cubic}), 
one finds that the $9$ remaining inequalities
independent of $\mu$ are satisfied by the 
following choice of coefficients:
\begin{align}
\bar{\alpha}_{0,0} & = 0, \; \; \;  
\bar{\alpha}_{1,0} = - \frac{5}{2}, \; \; \; 
\bar{\alpha}_{0,1} = - \frac{3}{2}, \; \; \; 
\bar{\beta}_{0,0} = 0, \; \; \; 
 \bar{\beta}_{1,0} = - 4, \; \; \; 
 \bar{\beta}_{0,1} = - \frac{5}{2}, \nonumber \\
\bar{\gamma}_{0,0} & = -2, \; \; \;  
\bar{\gamma}_{1,0} = 6, \; \; \;  
\bar{\gamma}_{0,1} = 1, \; \; \; 
\bar{\gamma}_{2,0} = 2, \; \; \;  
\bar{\gamma}_{1,1} = 0, \; \; \; 
\bar{\gamma}_{0,2} = - \frac{3}{4}. \nonumber 
\end{align}
The remaining $3$ inequalities from~(\ref{eq:R_1_cubic}),
which depend on $\mu$, are then given by 
\begin{align}
\frac{3}{2} - 32 \mu & > 0, \; \; \; 
\frac{17}{4} - 40 \mu > 0, \; \; \;
\frac{5}{2} + 80 \mu > 0,
\end{align}
and are satisfied for every $0 < \mu < 3/64$. 
Under this parameter choice, (\ref{eq:Escher_perturbed})
coincides with~(\ref{eq:cubic_R}), and is 
a cross-stencil reversible CDS 
by Lemma~\ref{lemma:reversible_1}. Furthermore, since the original
CDS~(\ref{eq:Escher}) has a stable hyperbolic algebraic limit cycle $\Gamma$, 
which is a regular oval of $\{h = 0\}$ (see~\cite{Escher} or 
Lemma~\ref{lemma:Escher}), Theorem~\ref{theorem:perturbations} 
guarantees that $\Gamma$ remains a hyperbolic algebraic limit cycle
with the same stability for every sufficiently small $\mu > 0$. 
\end{proof}

\noindent
\textbf{Remark}. A more detailed analysis shows that 
Theorem~\ref{theorem:cubic_R} holds for every $0 < \mu < 3/64$. 

Theorem~\ref{theorem:cubic_R} gives the following corollary.

\begin{corollary}
\label{corollary:cubic_R}
$H_{a,r}^c(3) \ge 1$. 
\end{corollary}

\subsection{Quartic system with two algebraic limit cycles}
By Lemmas~\ref{lemma:CDS_R_1} and~\ref{lemma:bHarnack_3}, 
there exists a reversible quartic CDS with an algebraic limit cycle.
We now use an analogous approach as in 
Section~\ref{sec:cancel_1} to construct 
a reversible system with two algebraic limit cycles. 

\begin{lemma}
\label{lemma:quartic_double_R}
Let $\Gamma_1$ and $\Gamma_2$ be the two ovals of $\{H = 0\}$, 
enclosing $(-1,0)$ and $(1,0)$, respectively, where
$H$ is the quartic polynomial given by~$(\ref{eq:H_double})$.
Consider the quartic \emph{DS}
\begin{align}
\frac{\mathrm{d} u}{\mathrm{d} t} 
& = \bar{\alpha}_{0,0} H(u,v) - 8 c(u,v) v, \nonumber \\
\frac{\mathrm{d} v}{\mathrm{d} t} 
& = 2 d(u,v) u^3 + g_1(u,v) (1 - 4 u^2 + 4 v^2) - 8 c(u,v) u,
\label{eq:quartic_double_R}
\end{align} 
where, for $\bar{\alpha}_{0,0}, \bar{\beta}_{i,j}, \bar{\gamma}_{i,j} \in \mathbb{R}$,  
\begin{align}
g_1(u,v) & = \bar{\beta}_{0,0} + \bar{\beta}_{1,0} u + \bar{\beta}_{0,1} v
+ \bar{\beta}_{2,0} u^2 + \bar{\beta}_{1,1} u v + \bar{\beta}_{0,2} v^2, \nonumber \\ 
c(u,v) & = \frac{1}{4} 
\left( d(u,v) -  g_1(u,v) u \right), \; \; \; 
d(u,v) = \bar{\gamma}_{0,0} + \bar{\gamma}_{1,0} u + \bar{\gamma}_{0,1} v.
\label{eq:quartic_double_fgc_R}
\end{align}
Define 
\begin{align}
F = \frac{1}{4}
\left(d \frac{\partial g_1}{\partial v} - g_1 \frac{\partial d}{\partial v}
+ \bar{\alpha}_{0,0} \left( g_1 + u \frac{\partial g_1}{\partial u} 
- \frac{\partial d}{\partial u}\right) \right).
\label{eq:F_R}
\end{align}
If $c \ne 0$ and $F \ne 0$ on 
$\overline{\Omega(\Gamma_1)}$ and $\overline{\Omega(\Gamma_2)}$,
then $\Gamma_1$ and $\Gamma_2$ are hyperbolic algebraic limit cycles 
of~$(\ref{eq:quartic_double_R})$. Furthermore, for each $i = 1,2$, 
the cycle $\Gamma_i$ is stable if $c F < 0$ on $\Gamma_i$, and 
unstable if $c F > 0$ on $\Gamma_i$.
\end{lemma}

\begin{proof}
Substituting~(\ref{eq:H_double}) into~(\ref{eq:DS_1})
and choosing $f_1 = \bar{\alpha}_{0,0}$, one obtains
\begin{align}
\frac{\mathrm{d} u}{\mathrm{d} t} 
& = \bar{\alpha}_{0,0} H(u,v) - 8 c(u,v) v, \nonumber \\
\frac{\mathrm{d} v}{\mathrm{d} t} 
& = 2 u^3 \left(g_1(u,v) u + 4 c(u,v) \right) 
+ g_1(u,v) (1 - 4 u^2 + 4 v^2) - 8 c(u,v) u.
\label{eq:quartic_double_R_pre}
\end{align} 
The first equation is at most quartic if $c$ is at most cubic. 
The second equation is then at most quartic if $g_1$ is at most quadratic 
and $g_1 u + 4 c$ is at most affine. Choosing~(\ref{eq:quartic_double_fgc_R}),
DS~(\ref{eq:quartic_double_R_pre}) becomes~(\ref{eq:quartic_double_R}).
Furthermore, function $F$, defined in~(\ref{eq:F}), becomes~(\ref{eq:F_R}).

Assume that $c \ne 0$ and $F \ne 0$ on 
$\overline{\Omega(\Gamma_1)}$ and $\overline{\Omega(\Gamma_2)}$.
Since $c$ and $F$ are continuous and non-zero, 
each function has a constant sign on each region. 
For each $i = 1,2$, choose a sufficiently small simply connected open set
$\mathbb{S}_i$ containing $\overline{\Omega(\Gamma_i)}$ on which $c \ne 0$
and $F \ne 0$, and which is disjoint from the other oval. 
By Theorem~\ref{theorem:AL_strong}, $\Gamma_1$ and $\Gamma_2$ are hyperbolic 
algebraic limit cycles of~(\ref{eq:quartic_double_R}).
Furthermore, $\{H = 0\}$ consists purely of two disjoint, non-nested regular ovals, 
and $H(-1,0) = H(1,0) = -1$. Hence $H < 0$ in $\Omega(\Gamma_i)$
and $H > 0$ immediately outside $\Gamma_i$. Therefore, 
$s(\mathbf{x}_{\Gamma_i}) = \textrm{sign}(c(\mathbf{p}_i))$ 
for any $\mathbf{p}_i \in \Gamma_i$. Using~(\ref{eq:K_F}), 
$\textrm{sign}(K_{\Gamma_i}) = \textrm{sign}(c(\mathbf{p}_i) F(\mathbf{p}_i))$,
giving the stated stability. 
\end{proof}

The coefficients in the following theorem 
can be obtained by solving the linear inequalities~(\ref{eq:reversible_1}),
as in the proof of Theorem~\ref{theorem:cubic_R}.
The phase-plane diagram of the constructed reversible 
quartic CDS~(\ref{eq:quartic_2_R}) is shown in Figure~\ref{fig:1}(d).

\begin{theorem}
\label{theorem:quartic_2_R}
System
\begin{align}
\frac{\mathrm{d} x}{\mathrm{d} t}
& = 73648 + 31653820 x + 76492 y
+ 10682202 x^2 - 18853513 x y + 2228656 y^2  \nonumber \\
& - 522256 x^3 + 359200 x^2 y + 480492 x y^2 + 2036 y^3
- 664 x^4 - 1204 x^3 y - 398 x^2 y^2 - 1167 x y^3, \nonumber \\
\frac{\mathrm{d} y}{\mathrm{d} t}
& = 28141104 + 22676288 x + 3748360 y
+ 23171387 x^2 - 33870534 x y + 5060064 y^2 \nonumber \\
& + 4096 x^3 - 1513312 x^2 y + 1868440 x y^2 - 33304 y^3
- 64 x^3 y - 1404 x^2 y^2 - 2398 x y^3 - 192 y^4,
\label{eq:quartic_2_R}
\end{align}
is a cross-stencil reversible \emph{CDS} of degree $4$. 
Furthermore, $\Gamma_1, \Gamma_2 \subset \{h = 0\} \cap \mathbb{R}_{>}^2$ 
are stable hyperbolic algebraic limit cycles
enclosing $(16,15)$ and $(8,9)$, respectively,
where $h(x,y) = H((12 - 2 x + y)/5, (-12 - 3 x + 4 y)/5)$ 
and $H$ is the quartic polynomial given by~$(\ref{eq:H_double})$.
\end{theorem}

\begin{proof}
Consider DS~(\ref{eq:quartic_double_R}) with 
$\bar{\alpha}_{0,0} = 17$,
\begin{align}
g_1(u,v) & = 1557 - 453 u + 104 v - 19 u^2 + 7 u v - \frac{7}{8} v^2, \nonumber \\ 
d(u,v) & = 368 - \frac{385}{16} u + \frac{2673}{64} v,
\end{align}
and $c(u,v) = (d(u,v) - u g_1(u,v))/4$. Computing $F$ from~(\ref{eq:F_R}), one obtains
\begin{align}
F(u,v) & = \frac{1}{2048} \left(61944 + 1838760 u + 575488 v 
- 176072 u^2 + 143416 u v - 26327 v^2 \right).
\end{align}
Since $H \le 0$ on $\overline{\Omega(\Gamma_1)}$ and $\overline{\Omega(\Gamma_2)}$, 
these regions satisfy $2 (u^2 - 1)^2 + 4 v^2 \le 1$, and hence 
$\sqrt{1 - 1/\sqrt{2}} \le |u| \le \sqrt{1 + 1/\sqrt{2}} < 4/3$ and $|v| \le 1/2$.
On $\overline{\Omega(\Gamma_1)}$, $-4/3 < u < -1/2$, $|v| \le 1/2$, and hence 
\begin{align}
2048 F(u,v) & <  - \frac{1422244}{3} < 0, \; \; \; \; 
d(u,v) > 368 - \frac{2673}{128} > 0, \; \; \; \; 
g_1(u,v) > 1557 - 52 - \frac{304}{9} - \frac{14}{3} - \frac{7}{32} > 0.
\end{align}
Hence, $c > 0$, so $c F < 0$ on $\overline{\Omega(\Gamma_1)}$.
Analogous argument shows that $c < 0$ and $F > 0$, and hence $c F < 0$,
on $\overline{\Omega(\Gamma_2)}$. Consequently, 
by Lemma~\ref{lemma:quartic_double_R}, $\Gamma_1$ and $\Gamma_2$ 
are stable hyperbolic algebraic limit cycles. 

Consider the affine change of coordinates
\begin{align}
x & = 12 - 4 u + v, \; \; \; \; y = 12 - 3 u + 2 v.
\end{align}
Since it is invertible, this map preserves the quartic degree,
the two hyperbolic algebraic limit cycles and their stability. 
Furthermore, since on both ovals $|u| < 4/3$ and $|v| \le 1/2$,
$x > 37/6 > 0$ and $y > 7 > 0$, i.e. both ovals are located
in the positive quadrant. Moreover, using the inverse map 
$u = (12 - 2 x + y)/5$, $v = (-12 - 3 x + 4 y)/5$ and 
rescaling time via $t \to t/4000$, DS~(\ref{eq:quartic_double_R})
under the stated choice of parameters maps to~(\ref{eq:quartic_2_R}).
The coefficients of this DS satisfy the reversibility
conditions~(\ref{eq:reversible_1}). Therefore, 
by Lemma~\ref{lemma:reversible_1}, 
(\ref{eq:quartic_2_R}) is a cross-stencil reversible CDS.
\end{proof}

Theorem~\ref{theorem:quartic_2_R} gives the following corollary.

\begin{corollary}
\label{corollary:quartic_2_R}
$H_{a,r}^c(4) \ge 2$. 
\end{corollary}

\section{Discussion}
\label{sec:6}
In this paper, we have developed general constructions
for realizing algebraic ovals as hyperbolic algebraic limit cycles 
in both general and reversible chemical dynamical systems. 
Through this work, we have established connections
between the topology of real algebraic curves and oscillatory
dynamics in chemical reaction networks, and derived
improved lower bounds on the number of algebraic
limit cycles in chemical dynamical systems, as 
summarized in Table~\ref{tab:2}.
These results contribute to a broader framework for
understanding how chemical constraints shape the 
admissible dynamics, complementing previous results
on equilibria, oscillations and chaos in chemical dynamical systems~\cite{CEpaper,
Homoclinic,Test_models,QCM,Chaos_1,Chaos_2,Yeung, George_Jiaxin_Polly, George_Balazs_Polly, George_Pantea}. Our results also raise a number of open questions,
some of which we now discuss.

In Theorem~\ref{theorem:semi_chemical_Harnack_stronger} in Section~\ref{sec:4}, 
we have constructed semi-chemical polynomials of every degree $n \ge 3$ 
with a single leading monomial
and $M(n) - (1 + (-1)^n)/2$ ovals in $\mathbb{R}_{>}^2$, 
allowing us to prove in Theorem~\ref{theorem:semi_chemical_improved}
that $H_{a}^c(n) \ge M(n) - (1 + (-1)^n)/2$ for every $n \ge 5$.
In particular, for odd $n \ge 5$, the bound 
$H_{a}^c(n) \ge M(n)$ matches the best known
lower bound for $H_a(n)$; see Table~\ref{tab:2}.
However, for even $n \ge 6$, the established bound 
$H_{a}^c(n) \ge M(n) - 1$ includes one limit cycle less
than the best known lower bound $H_a(n) \ge M(n)$. 

\textbf{Problem $1$}. \emph{Prove or disprove: $H_{a}^c(n) \ge M(n)$ 
for every even $n \ge 6$.}

\noindent 
One way to address this problem is via chemical Harnack M-curves. 
In particular, in Section~\ref{sec:4}, we have shown that 
quadratic chemical Harnack M-curves do not exist in Lemma~\ref{lemma:cHarnack_2}, 
and constructed chemical Harnack M-curves of degrees three and four
in Lemmas~\ref{lemma:cHarnack_3} and~\ref{lemma:cHarnack_4}, respectively.
If such curves exist for every $n \ge 6$, then 
Lemma~\ref{lemma:chemical} would imply that $H_{a}^c(n) \ge M(n)$ 
for every even $n \ge 6$. 

\textbf{Problem $1^*$}. \emph{Prove or disprove: chemical 
Harnack M-curves exist in every even degree $n \ge 6$.}

In Lemma~\ref{lemma:buffered_construct} in Section~\ref{sec:5},
we have shown that every polynomial of degree $n$
can be transformed into a buffered polynomial of degree $n+2$
while preserving its regular ovals and placing them in the positive quadrant.
Theorem~\ref{theorem:buffered_n} then gives $H_{a,r}^c(n) \ge M(n-3)$.

\textbf{Problem $2$}. \emph{Determine the smallest $k \in \{0,1,2,3\}$ such that $H_{a,r}^c(n) \ge M(n-k)$ for every $n \ge 4$.}

One way to establish this inequality with $k = 1$ is via buffered Harnack M-curves.
In particular, in Section~\ref{sec:5}, we have shown that 
quadratic buffered Harnack M-curves do not exist 
in Lemma~\ref{lemma:bHarnack_2}, and constructed 
a cubic buffered Harnack M-curve in Lemma~\ref{lemma:bHarnack_3}.
If such curves exist for every $n \ge 4$, then 
Lemma~\ref{lemma:CDS_R_1} would imply that 
$H_{a,r}^c(n) \ge M(n-1)$ for every $n \ge 4$.

\textbf{Problem $2^*$}. \emph{Prove or disprove: buffered Harnack 
M-curves exist in every degree $n \ge 4$.}

In Section~\ref{sec:5}, we have considered only
those reversible CDSs that admit a CRN realization with a cross-stencil,
in which each vertex of the E-graph is connected to all its
nearest horizontal and vertical neighbors, 
and to no other vertices, within the triangular lattice.
For this subclass of reversible systems, we have derived 
the lower bounds on $H_{a,r}^c$ presented in Table~\ref{tab:2}.
A natural extension is to investigate reversible CDSs with other stencils, 
including those with diagonal edges or longer-range connections.

More broadly, it is of interest to investigate CDSs
subject to weaker constraints on their coefficients
arising from reaction-network connectivity.
In particular, two additional Hilbert numbers are 
defined in~\cite{CEpaper}: the supremum of the number of algebraic limit cycles 
in $\mathbb{R}_{>}^2$ among $n$-degree planar CDSs
admitting weakly reversible reaction-network realizations,
which we denote by $H_{a,w}^c(n)$, and the corresponding
supremum among systems admitting CRN realizations 
in which every reaction is of the form 
$i X + j Y \xrightarrow[]{} i' X + j' Y$
with $i + j \le n$ and $i' + j' \le n$, 
which we denote by $H_{a,m}^c(n)$.
These Hilbert numbers satisfy
$H_{a,r}^c(n) \le H_{a,w}^c(n) \le H_{a,m}^c(n) \le H_{a}^c(n)$
for every $n \ge 2$. In this paper, we have focused 
on the two extremes of this hierarchy. 
Determining sharper bounds on these chemical Hilbert numbers 
and characterizing their relationships remain open problems.

Finally, the focus of this paper has been 
on the number of chemical algebraic limit cycles. 
Nevertheless, most of our results also address their stability
and uniqueness, and we have sought to construct examples
with larger numbers of stable hyperbolic limit cycles.
Characterizing the possible arrangements 
and stability properties of algebraic limit cycles in chemical
dynamical systems remains an open problem.

\vspace{0.5cm}
\noindent
\textbf{Funding}. Tomislav Plesa would like to thank Peterhouse, University of Cambridge, for a fellowship. He also acknowledges partial support from Research England through the Expanding Excellence in England (E3) funding stream awarded to MARS (Mathematics for AI in Real-world Systems) in the School of Mathematical Sciences at Lancaster University. This work was partially supported by the Ministerio de Ciencia e Innovaci\'{o}n, grant PID2025-168432NB-I00B, awarded to Armengol Gasull. Gheorghe Craciun was supported in part by the US National Science Foundation grant DMS-2051568, and by the Simons Foundation Travel Grant SFI-MPS-TSM-00025952.

\vspace{0.5cm}
\noindent
\textbf{AI use declaration}. During the preparation of this manuscript, 
which has extended over more than a year, the authors used AI tools 
to assist in checking and verifying mathematical results. The authors
take full responsibility for the content of the manuscript.

\end{document}